\documentclass[10pt]{amsart}
\usepackage{amsfonts,amsmath,amsthm,amssymb,color}
\usepackage[all]{xy}
\usepackage[pagebackref, colorlinks=true, linkcolor=blue, citecolor=blue]{hyperref}

\numberwithin{equation}{section}

\DeclareMathOperator*{\colim}{colim}

\newtheorem{theorem}{Theorem}[section]
\newtheorem{lemma}[theorem]{Lemma}

\newtheorem{proposition}[theorem]{Proposition}

\newtheorem{corollary}[theorem]{Corollary}

\theoremstyle{definition}
\newtheorem{definition}[theorem]{Definition}
\newtheorem{example}[theorem]{Example}

\theoremstyle{remark}
\newtheorem{remark}[theorem]{Remark}

\newcommand{\id}{\operatorname{id}}

\newcommand{\Mor}{\operatorname{Mor}}
\newcommand{\Ext}{\operatorname{Ext}}
\newcommand{\Irr}{\operatorname{Irr}}
\newcommand{\coker}{\operatorname{coker}}
\newcommand{\Lt}{\operatorname{L}}
\newcommand{\Mt}{\operatorname{M}}
\newcommand{\Ob}{\operatorname{Ob}}
\newcommand{\Paths}{\operatorname{Paths}}

\newcommand{\GProj}{\operatorname{\mathbf{GProj}}}
\newcommand{\Gproj}{\operatorname{\mathbf{Gproj}}}
\newcommand{\GInj}{\operatorname{\mathbf{GInj}}}
\newcommand{\projdim}{\operatorname{projdim}}
\newcommand{\injdim}{\operatorname{injdim}}
\newcommand{\Mod}{\operatorname{\mathbf{Mod}}}
\newcommand{\fgMod}{\operatorname{\mathbf{mod}}}
\newcommand{\Ho}{\operatorname{\mathbf{Ho}}}
\newcommand{\Smod}{\operatorname{\mathbf{\underline{Mod}}}}

\newcommand{\Hom}{\operatorname{Hom}}

\newcommand{\Rep}{\operatorname{\mathbf{Rep}}}
\newcommand{\rep}{\operatorname{\mathbf{rep}}}

\begin{document}

\title{Quiver representations and Gorenstein-projective modules}
\author{Francesco Meazzini}

 \address{\newline
Universit\`a degli studi di Roma La Sapienza,\hfill\newline
Dipartimento di Matematica \lq\lq Guido
Castelnuovo\rq\rq,\hfill\newline
P.le Aldo Moro 5,
I-00185 Roma, Italy.}
\email{meazzini@mat.uniroma1.it}
\urladdr{www.mat.uniroma1.it/people/meazzini/}

\renewcommand{\subjclassname}{%
\textup{2010} Mathematics Subject Classification}


\date{July 1, 2020.}

\begin{abstract}
We consider a finite acyclic quiver $\mathcal{Q}$ and a quasi-Frobenius ring $R$. We endow the category of quiver representations over $R$ with a model structure, whose homotopy category is equivalent to the stable category of Gorenstein-projective modules over the path algebra $R\mathcal{Q}$.

As an application, we then characterize Gorenstein-projective $R\mathcal{Q}$-modules in terms of the corresponding quiver $R$-representations; this generalizes a result obtained by Luo-Zhang to the case of not necessarily finitely generated $R\mathcal{Q}$-modules, and partially recover results due to Enochs-Estrada-Garc\'ia Rozas, and to Eshraghi-Hafezi-Salarian. Our approach to the problem is completely different since the proofs mainly rely on model category theory.
\end{abstract}

\maketitle
\tableofcontents

\section*{Introduction}

The notion of finitely generated Gorenstein-projective modules dates back to 1966. Following Auslander~\cite{Aus67}, a finitely generated module $M$ over a commutative and Noetherian ring $k$ is called Gorenstein-projective if it satisfies the conditions below:
\begin{itemize}
\item $\Ext^i_k(M, k) = \Ext^i_k(\Hom_k(M, k), k) = 0$ for all $i>0$,
\item the natural biduality homomorphism $M \to \Hom_k(\Hom_k(M, k), k)$ is an isomorphism.
\end{itemize}

We will denote by $\Gproj(k)$ the full subcategory of left $k$-modules $\Mod(k)$, whose objects are finitely generated Gorenstein-projective modules. This notion extends the one of finitely generated projective ones, and led Auslander to introduce the notion of $G$-dimension for any finitely generated module $M$, which is defined to be the minimal length of a ``Gorenstein-projective resolution" of $M$, see \cite{Aus67}. Subsequently Auslander and Bridger developed the theory only assuming the ring $k$ to be associative and both left and right Noetherian, see \cite{AB69}.

Few years later Iwanaga introduced the general notion of Gorenstein ring (see \cite{Iw79} and \cite{Iw80}) generalizing the definition of the commutative case. Some authors refer to non-commutative Gorenstein rings as \emph{Iwanaga-Gorenstein} rings. Examples of Gorenstein rings are quasi-Frobenius rings and group rings $k[G]$ for any commutative Gorenstein ring $k$ and any finite group $G$, see \cite{EN55}.

Auslander's notion of finitely generated Gorenstein-projective modules still plays an interesting role in both Algebra and Geometry. For instance, a famous result due to Buchweitz~\cite{Bu87} states an equivalence of categories
\[ D_{sing}^{b}(k) \simeq \underline{\Gproj}(k) \]
where:
\begin{description}
\item[1] $k$ is a (not necessarily commutative) Gorenstein ring,
\item[2] $D_{sing}^{b}(k)$ denotes the category of singularities, i.e. the Verdier quotient of the bounded derived category $D^b(\fgMod(k))$ modulo the subcategory of perfect complexes, and $\fgMod(k)$ denotes the category of finitely generated $k$-modules,
\item[3] $\underline{\Gproj}(k)$ denotes the stable category of finitely generated Gorenstein-projective $k$-modules, i.e. its objects are the same as $\Gproj(k)$, while the morphisms are defined as
\[ \Hom_{\underline{\Gproj}(k)}(M, N) = \frac{\Hom_{\Gproj(k)}(M, N)}{\left\{ f\colon M\to N \,\,\vert\,\, f \mbox{ factors trough a projective $k$-module}\right\}}. \]
\end{description}
Nevertheless, at that time the so-called Gorenstein homological algebra still presented the problem of being generalized to not necessarily finitely generated modules.

In 1995, in \cite{EJ95}, Enochs and Jenda defined Gorenstein-projective and Gorenstein-injective modules over an arbitrary associative ring $k$.

\begin{definition}\label{EJ}
Let $k$ be an associative ring. A $k$-module $M\in\Mod(k)$ is called \emph{Gorenstein-projective} if there exists an exact sequence of projective modules
\[ \cdots\to P^{-1}\xrightarrow{d^{-1}} P^0\xrightarrow{d^0} P^1\to \cdots \]
that remains exact under the functor $\Hom(-,P)$ for every projective $k$-module $P\in\Mod(k)$, and such that $M\cong \ker\{d^0\}$.

Dually, $M\in\Mod(k)$ is called \emph{Gorenstein-injective} if there exists an exact sequence of injective modules
\[ \cdots\to J^{-1}\xrightarrow{d^{-1}} J^0\xrightarrow{d^0} J^1\to \cdots \]
that remains exact under the functor $\Hom(J,-)$ for every injective $k$-module $J\in\Mod(k)$, and such that $M\cong \ker\{d^0\}$.
\end{definition}

Finally Avramov, Buchweitz, Martsinkovsky, and Reiten proved that over an associative ring $k$ which is both left and right Noetherian, for any finitely generated $k$-module, Definition~\ref{EJ} and Auslander's notion of Gorenstein-projectivity coincide, see \cite{ABMR}. To the author knowledge \cite{ABMR} is still not published. With regard to the result mentioned above we then refer to \cite{Ch00}, where it was included (see Theorem 4.2.6) by Christensen for a commutative ring $k$, but it is straightforward to check that the argument works also in the general case.

The basic properties of Gorenstein-projective and Gorenstein-injective modules were then investigated by Enochs and Jenda. We refer to \cite{EJ11} for a complete exposition on the subject.

In general, the problem of explicitly describing Gorenstein-projective $k$-modules presents several difficulties, but for particular choices of the base ring $k$ the situation becomes extremely simple. For instance, if $k$ is a quasi-Frobenius ring (such as $k=\mathbb{C}[\epsilon]$, the algebra of dual numbers) then every $k$-module is Gorenstein-projective, namely $\Mod(k)=\GProj(k)$. During the last decade, authors tried to avoid problems related to the base ring, trying to increase their understanding of the class of Gorenstein-projective $k$-modules through the information they had about the same class over a simpler ring $R$. The key in this approach relies on the equivalence of categories
\[ \Rep(\mathcal{Q}, R) \simeq \Mod(\Lambda^{op}) \]
where $\mathcal{Q}$ is a finite acyclic quiver, $\Lambda^{op}$ is the (opposite) path algebra of $\mathcal{Q}$ over a unitary ring $R$, and $\Rep(\mathcal{Q}, R)$ denotes the category of representations of $\mathcal{Q}$ over $R$. It is worth noticing that even if $R=\mathbb{C}[\epsilon]$ is the self-injective algebra of dual numbers, the path algebra $\Lambda^{op}$ is a non-commutative $1$-Gorenstein ring (see Lemma~\ref{1gor} below). Hence, the hope is to find an explicit description for the subcategory $\GProj(\Lambda^{op})\subseteq \Mod(\Lambda^{op})$ in terms of quiver representations with values in the easier class $\GProj(\mathbb{C}[\epsilon])=\Mod(\mathbb{C}[\epsilon])$ through the equivalence of categories above. This approach has been followed in a very general setting by Luo and Zhang in \cite{LZ13}. More precisely, they restrict their attention to the class of finitely generated Gorenstein-projective $\Lambda^{op}$-modules, but on the other hand their result allows $R$ to be any finite-dimensional algebra over a field. As an application we extend their result to the whole category of Gorenstein-projective $\Lambda^{op}$-modules, whenever $R$ is a quasi-Frobenius ring, see Corollary~\ref{gorcof}.
It is important to point out that Gorenstein-projective modules have already been described by Enochs, Estrada and Garc\'ia Rozas,~\cite[Corollary 6.4]{EEG09}, if the ring is Gorenstein; later Eshraghi, Hafezi and Salarian,~\cite{EHS13}, obtained the same characterization for any ring.
The main innovation of this paper concerns the different approach to the problem, which does not involve standard techniques of quiver representations, but relies instead on model category theory.

The reason for our assumption on $R$ is not technical, meaning that up to mild arrangments our argument works in a more general setting (namely when $R$ is a possibly non-commutative Gorenstein ring, see Remark~\ref{generalization}). Nevertheless, if $R$ is quasi-Frobenius, the class $\mathcal{W}$ of Reedy weak equivalences in $\Rep(\mathcal{Q}, R)$ admits an easy characterization (see Theorem~\ref{modelrep}), and the stable category $\underline{\GProj}(\Lambda^{op})$ is equivalent to the localization $\Rep(\mathcal{Q}, R)[\mathcal{W}^{-1}]$. Moreover, the $RP$-cofibrant $R$-representations can be more easily described, see Theorem~\ref{modelrep}.
In particular, when $R=\mathbb{C}[\epsilon]$ is the algebra of dual numbers, the stable category $\underline{\Gproj}(\Lambda^{op})$ of finitely generated Gorenstein-projective $\Lambda^{op}$-modules has been deeply investigated in \cite{RZ13}.

The plan of the paper is as follows. We begin Section~\ref{tmcoqroaqfr} by recalling the stable model structure of modules over a quasi-Frobenius ring. Here we only assume the reader to be familiar with the basic notions of model category theory, for which we refer to~\cite{Hir03,Hov99}. We proceed by briefly surveying the basic definitions concerning quiver representations in order to obtain our first result, which states the existence of a model structure on the category $\Rep(\mathcal{Q}, R)$ of quiver representations over a quasi-Frobenius ring $R$.

\begin{theorem}[see Theorem~\ref{modelrep}]
Let $\mathcal{Q}$ be a finite acyclic quiver, and let $R$ be a quasi-Frobenius ring. Then the category $\Rep(\mathcal{Q}, R)$ admits a model structure, called the \emph{Reedy projective} model structure, where a morphism between $R$-representations $M\to N$ is:
\begin{description}
\item[1] a weak equivalence (called Reedy stable equivalence) if it is a pointwise stable equivalence in $\Mod(R)$, i.e. if for every vertex $j\in \mathcal{Q}_0$ the morphism $M_j\to N_j$ is a stable equivalence of $R$-modules, 
\item[2] a fibration (called $RP$-fibration) if it is a pointwise surjection in $\Mod(R)$, i.e. if for every vertex $j\in \mathcal{Q}_0$ the morphism $M_j\to N_j$ is surjective, 
\end{description}
Moreover, an $R$-representation $M$ is $RP$-cofibrant if and only if for every vertex $j\in \mathcal{Q}_0$ the natural morphism of $R$-modules
\[ \bigoplus\limits_{\stackrel{\alpha\in \mathcal{Q}_1}{\tau(\alpha) = j}}M_{\sigma(\alpha)} \to M_j \]
is injective.
\end{theorem}

The importance of the above result lies in the explicit characterization of the $RP$-cofibrant representations; these will correspond to Gorenstein-projective $\Lambda^{op}$-modules. 

In Theorem~\ref{modelrep} we also give a characterization of cofibrations in the model structure above, and we describe a dual model structure, which will be called the \emph{Reedy injective} model structure. Except for the explicit description of cofibrant objects, the result above almost immediately follows from the standard Reedy model structures on categories of diagrams (see Theorem~\ref{Reedy}), whence we decided to preserve the name. This result plays a crucial role in what follows, since the cofibrant $R$-representations corresponds to Gorenstein-projective modules over the (opposite) path algebra of $\mathcal{Q}$ over $R$, giving an explicit description of such modules in terms of $\mathcal{Q}$.
It should be noticed that finitely generated cofibrant $R$-representations are precisely the so-called \emph{monic representations} introduced by Luo and Zhang in \cite{LZ13}.

The purpose of Section~\ref{searse} is to find a different description of stable equivalences and Reedy stable equivalences, in order to prove the main results in the following sections. We also investigate the relation between projective $R$-representations and cofibrant ones. We prove that projective $R$-representations are just cofibrant $R$-representations which are vertexwise projective over $R$.

\begin{lemma}[see Lemma~\ref{projrep}]
Let $\mathcal{Q}$ be a finite acyclic quiver, and let $R$ be a quasi-Frobenius ring. Then the following are equivalent for an $R$-representation $M$.
\begin{description}
\item[1] $M$ is a projective object in $\Rep(\mathcal{Q},R)$,
\item[2] $M$ is $RP$-cofibrant and $M_j$ is a projective $R$ module for every $j\in\mathcal{Q}_0$.
\end{description}
\end{lemma}

Lemma~\ref{projrep} also gives an equivalent description for injective $R$-representations in terms of $RI$-fibrant ones.

In Section~\ref{gorenstein} we present our main results. First we recall an important result due to Hovey, \cite{Hov02}, where the category of modules over a (not necessarily commutative) Gorenstein ring $G$ is endowed with a model structure in which cofibrant objects are precisely the Gorenstein-projective $G$-modules. We proceed by showing that, under mild assumptions, the (opposite) path algebra $\Lambda^{op}$ of a quiver $\mathcal{Q}$ over a quasi-Frobenius ring $R$ is $1$-Gorenstein (see Lemma~\ref{1gor}). In particular, this permits to transfer Hovey's model structure on the category of $R$-representations of $\mathcal{Q}$ through the well known equivalence $\Rep(\mathcal{Q}, R)\simeq\Mod(\Lambda^{op})$. We then prove our main theoretic result, which will lead immediately to the main application of the paper, see Corollary~\ref{lz}.

\begin{theorem}[see Theorem~\ref{modelcoincide}]
Let $\mathcal{Q}$ be a finite acyclic quiver and let $R$ be a quasi-Frobenius ring. Then the (transferred) Hovey-projective model structure on $\Rep(\mathcal{Q},R)$ and the Reedy-projective model structure coincide.
\end{theorem}

The reader may notice that Theorem~\ref{modelcoincide} presents a dual statement about the Hovey-injective model structure and the Reedy-injective one. As an immediate consequence of the theorem above, we obtain that the cofibrant objects with respect to the Hovey model structure (i.e. the Gorenstein-projective $\Lambda^{op}$-modules) correspond to the $RP$-cofibrant representations, as explicitly stated below.

\begin{corollary}[see Corollary~\ref{gorcof}]\label{lz}
Let $\mathcal{Q}$ be a finite acyclic quiver and let $R$ be a quasi-Frobenius ring. Consider the path algebra $\Lambda = R\mathcal{Q}$. A module $M\in\Mod(\Lambda^{op})$ is Gorenstein-projective if and only if the corresponding $R$-representation  $M\in\Rep(\mathcal{Q},R)$ satisfies the following condition:

the morphism of $R$-modules
\[ \bigoplus\limits_{\stackrel{\alpha\in \mathcal{Q}_1}{\tau(\alpha) = j}}M_{\sigma(\alpha)} \to M_j \]
is injective for every vertex $j\in \mathcal{Q}_0$.
\end{corollary}

Again, Corollary~\ref{gorcof} states a dual result, characterizing Gorenstein-injective $\Lambda^{op}$-modules as $RI$-fibrant $R$-representations.
The equivalence $\Rep(\mathcal{Q}, R)\simeq\Mod(\Lambda^{op})$ restricts to an equivalence of categories
$\rep(\mathcal{Q}, R)\simeq\fgMod(\Lambda^{op})$ between finitely generated representations over $R$ and finitely generated $\Lambda^{op}$-modules. In particular, we recover the characterization already obtained in \cite{EEG09,EHS13,LZ13}; notice that the proofs of the cited papers are based on a completely different point of view.

In Section~\ref{homotopy} we investigate the main properties of the homotopy category $\Ho(\Rep(\mathcal{Q}, R))$, in particular we give an elementary proof that there exist equivalences of triangulated categories
\[ \underline{\GProj}(\Lambda^{op})\simeq\Ho(\Rep(\mathcal{Q}, R))\simeq\underline{\GInj}(\Lambda^{op}) \]
where $\underline{\GProj}(\Lambda^{op})$ denotes the stable category of Gorenstein-projective $\Lambda^{op}$-modules while, respectively, $\underline{\GInj}(\Lambda^{op})$ denotes the stable category of Gorenstein-injective ones.

\section{Model structures on quiver representations over a quasi-Frobenius ring}\label{tmcoqroaqfr}

This section is devoted to the description of two model structures on the category of quiver representations over a quasi-Frobenius ring.
As we will see, this yields two (nontrivial) model structures on the category of modules over a large class of (not necessarily quasi-Frobenius) rings, see Theorem~\ref{modelmodule}. If not specified, modules are assumed to be \emph{left} modules. The category of left $R$-modules will be denoted by $\Mod(R)$ for any unitary (not necessarely commutative) ring $R$. Then the category $\Mod(R^{op})$ is the category of right $R$-modules. We stress the fact that we do not restrict our interest to finitely generated modules.
We begin by recalling the standard model structure on the category of modules over a quasi-Frobenius ring, which is described in \cite{Hov99}, and in a more general setting in \cite{Pir86}.

\begin{definition}
A Noetherian (not necessarely commutative) ring $R$ is \emph{quasi-Frobenius} if it is injective both as a left and right $R$-module.
\end{definition}

In \cite{Fa66,FW67} Faith and Walker proved that the following conditions are equivalent:
\begin{description}
\item[1] $R$ is quasi-Frobenius,
\item[2] each projective right $R$-module is injective,
\item[2*] each injective right $R$-module is projective,
\item[3] each injective left $R$-module is projective,
\item[3*] each projective left $R$-module is injective.
\end{description}
It follows that a Noetherian ring $R$ is \emph{quasi-Frobenius} if and only if the classes of projective and injective $R$-modules coincide.

One of the more interesting classes of quasi-Frobenius rings is the one of self-injective algebras over a field $\mathbb{K}$.
Recall that a commutative $\mathbb{K}$-algebra $R$ is called \emph{self-injective} if it is an injective $R$-module.
Examples of self-injective $\mathbb{C}$-algebras are $\frac{\mathbb{C}[t]}{(t^n)}$, $n\in\mathbb{N}$.

Recall that given morphisms $f, g\colon M\to N$ of $R$-modules, then $f$ is called \emph{stably equivalent} to $g$ if the map $(f-g)$ factors through a projective $R$-module.

\begin{definition}\label{stable}
Let $R$ be a ring. The \emph{stable category} of $R$-modules is the category $\Smod(R)$ whose objects are left $R$-modules and whose morphisms are stable equivalence classes of morphisms in $\Mod(R)$. We shall call a morphism in $\Mod(R)$ a \emph{stable equivalence} if its class represents an isomorphism in $\Smod(R)$.
\end{definition}

\begin{remark}
We clearly have a functor $\gamma\colon \Mod(R)\to\Smod(R)$ that is the identity on objects. The category $\Smod(R)$ satisfies the following universal property. Given a category $\mathbf{C}$ and a functor $F\colon \Mod(R)\to \mathbf{C}$ such that for every projective/injective $R$-module $M$ there exists an isomorphism $F(M)\xrightarrow{\cong} F(0)$ in $\mathbf{C}$, then there exists a unique functor $G\colon \Smod(R)\to\mathbf{C}$ such that $F=G\circ \gamma$.
\end{remark}

\begin{theorem}\cite[Theorem 2.2.12]{Hov99}\label{Mark}
Suppose $R$ is a quasi-Frobenius ring. Then there is a model structure on $\Mod(R)$, where the cofibrations are injections, the fibrations are surjections, and the weak equivalences are the stable equivalences.
\end{theorem}

Notice that in the model structure of Theorem~\ref{Mark} every $R$-module is both fibrant and cofibrant, i.e. the morphism $0\to M$ is a cofibration and the morphism $M\to 0$ is a fibration for every module $M\in\Mod(R)$.

We now briefly recall the basic notions concerning quiver representations.

%\begin{definition}\label{kronecker}
%The Kronecker category $\mathbf{K}$ is defined as follows.
%\begin{description}
%\item[1] $\Ob(\mathbf{K}) = \{0,1\}$.
%\item[2] There are only two non-identity morphisms in $\mathbf{K}$, both of them with $1$ as domain and $0$ as codomain. 
%\end{description}
%\end{definition}

A \emph{quiver} can be simply understood as an oriented graph;
%and it may be defined as a functor $\mathcal{Q}\colon\mathbf{K}\to\mathbf{Sets}$. Notice that such a functor is equivalent 
more precisely a quiver $\mathcal{Q}$ is the data of two sets $\mathcal{Q}_0$ and $\mathcal{Q}_1$, together with two maps $\sigma,\tau\colon \mathcal{Q}_1\to \mathcal{Q}_0$.
Given a quiver $\mathcal{Q}$ we shall call $\mathcal{Q}_1$ the set of \emph{arrows}, and $\mathcal{Q}_0$ the set of \emph{vertices}. We will denote by $\sigma,\tau\colon \mathcal{Q}_1\to \mathcal{Q}_0$ the \emph{source} and \emph{target} map respectively. Therefore an arrow is an element $\alpha\in \mathcal{Q}$, and may be also denoted by $\sigma(\alpha) \xrightarrow{\alpha} \tau(\alpha)$. A quiver $\mathcal{Q}$ will be called \emph{acyclic} if it contains no closed paths, i.e. it does not exist a (non trivial) path $i_0\xrightarrow{f_1} i_1\to \cdots\to i_{n-1}\xrightarrow{f_n} i_n$ such that $i_0=i_n$, for any $n\geq 1$. In the literature, acyclic quivers are also called \emph{directed} (see e.g. \cite{RZ13}), while the name \emph{acyclic} can be found in \cite{LZ13}, for example. Since acyclic quivers will be dealing with direct and inverse categories, we will never use the name ``directed quiver". Example~\ref{fdrc} and Proposition~\ref{ffdrc} will clarify our choice.

\begin{definition}
Let $s\to t$ be a morphism in a small category $\mathbf{C}$. A factorization
\[ s \xrightarrow{f_1}\cdots\xrightarrow{f_n} t \]
is called a \emph{trivial factorization} if at least $n-1$ out of $\{f_i\}_{i\in\{1,\dots, n\}}$ are identity morphisms in $\mathbf{C}$.
A non-identity morphism in $\mathbf{C}$ is called \emph{irreducible} if it does not admit any (non trivial) factorization in $\mathbf{C}$.
\end{definition}

We shall remark that by definition identity morphisms are \emph{not} considered irreducible morphisms.
%For instance, in the Kronecker category $\mathbf{K}$ defined in \ref{kronecker} both non-identity morphisms are irreducible.

\begin{definition}
Let $s\to t$ be a non-identity morphism in a small category $\mathbf{C}$. A factorization
\[ s \xrightarrow{f_1}\cdots\xrightarrow{f_n} t \]
of $s\to t$ is called \emph{irreducible} if each morphism $f_i$ is \emph{irreducible} in $\mathbf{C}$.
\end{definition}

Recall that a category $\mathbf{C}$ is said to be \emph{finite} if $\Ob(\mathbf{C})$ is a finite set, and for every $s,t\in\mathbf{C}$ the set $\Hom_{\mathbf{C}}(s,t)$ is finite.
Clearly, given a finite category $\mathbf{C}$ every morphism admits at least one irreducible factorization. This naturally leads to the notion of \emph{free categories}, which we now introduce.

\begin{definition}
A small category $\mathbf{C}$ is said to be \emph{free} if every morphism in $\mathbf{C}$ admits a unique irreducible factorization.
\end{definition}

All the free categories that will be considered in the following turn out to be finite. A \emph{degree function} on a small category $\mathbf{C}$ is a map of sets $d\colon \Ob(\mathbf{C})\to \mathbb{N}$.

\begin{remark}
There exists a more general notion of degree functions, see \cite{Hir03}. Nevertheless we decided to present it as simply as possible, since it is not going to play a crucial role in what follows.
\end{remark}

\begin{definition}
A free \emph{direct} Reedy category $(\mathbf{R},d)$ consists of a free category $\mathbf{R}$ endowed with a degree function $d\colon\Ob(\mathbf{R})\to \mathbb{N}$ such that every non-identity morphism increases the degree, i.e. for every non-identity morphism $r_1\to r_2$ in $\mathbf{R}$ we have $d(r_1) < d(r_2)$.
Dually, a free \emph{inverse} Reedy category $(\mathbf{R},d)$ consists of a free category $\mathbf{R}$ endowed with a degree function $d\colon\Ob(\mathbf{R})\to \mathbb{N}$ such that every non-identity morphism decreases the degree.
\end{definition}

These are very special classes of Reedy categories, hence we decided to preserve the name. For the general notion of Reedy categories we refer to \cite{Hir03}.
We will be interested just in finite free direct (or inverse) Reedy categories.
Example~\ref{fdrc} and Proposition~\ref{ffdrc} will explain how these categories are related to quivers.
The adjective \emph{free} is due in fact to the correspondence with quivers.

\begin{example}\label{fdrc}
Let $\mathcal{Q}$ be a finite acyclic quiver. Then $\mathcal{Q}$ induces a finite free direct Reedy category $(\mathbf{Q}, d)$ with the following definitions.
\begin{description}
\item[1] $\Ob(\mathbf{Q}) = \mathcal{Q}_0$, i.e. the objects of $\mathbf{Q}$ are the vertices of the quiver. 
\item[2] $\Mor(\mathbf{Q}) = \Paths(\mathbf{Q})$, i.e. the morphisms in the category $\mathbf{Q}$ is the set of arrows $\mathcal{Q}_1$ plus all the possible compositions between them and the identity morphisms.
\item[3] A degree function can be defined as follows. Since the quiver $\mathcal{Q}$ is finite and acyclic, we have at least one vertex which is not the target of any arrow. We shall call such a vertex a \emph{source vertex}. Let us denote by $\{s_1,\dots ,s_n\}$ the set of source vertices in $\mathcal{Q}$. We then define $d(s_i)=0$ for every $i\leq n$. Now pick a vertex $j\in \mathcal{Q}_0$, and consider the set
\[ P_j\colon=\{s_i\to\cdots\to j\vert 1\leq i\leq n\} \]
of all paths starting from a source vertex and ending with the fixed vertex $j$. Since the quiver $\mathcal{Q}$ is acyclic and finite, $P_j$ is a finite (non-empty) set. Hence there exists at least one path of maximal length $(s_{i}\to j_1\to \cdots\to j_m= j)$ in $P_j$, for some $i\leq n$. We define $d(j)=m$.
\end{description}
\end{example}

From now on, given a finite acyclic quiver $\mathcal{Q}$ we will denote by $\mathbf{Q}$ the associated free direct Reedy category as explained in Example~\ref{fdrc}.

\begin{remark}
Let $\mathcal{Q}$ be a finite acyclic quiver. Thanks to Example~\ref{fdrc} we can associate to $\mathcal{Q}$ a finite free \emph{inverse} Reedy category $(\stackrel{\leftarrow}{\mathbf{Q}}, \stackrel{\leftarrow}{d})$ as follows. The category $\stackrel{\leftarrow}{\mathbf{Q}}$ is the same category as $\mathbf{Q}$, while the degree function is defined as
\[ \begin{aligned}
\stackrel{\leftarrow}{d}\colon \stackrel{\leftarrow}{\mathbf{Q}} &\to\mathbb{N}\\
j &\mapsto M - d(j)
\end{aligned} \]
where $d$ is the degree function defined in Example~\ref{fdrc} and
\[ M = \max\{d(j) \,\,\vert\,\, j\in Q_0\} \in\mathbb{N} \; . \]
\end{remark}

\begin{proposition}\label{ffdrc}
Every finite free direct Reedy category is (except for the degree function) of the form $\mathbf{Q}$ for some finite acyclic quiver $\mathcal{Q}$. Dually, every finite free inverse Reedy category is (except for the degree function) of the form $\stackrel{\leftarrow}{\mathbf{Q}}$ for some finite acyclic quiver $\mathcal{Q}$.
\begin{proof}
Given a finite free direct Reedy category $\mathbf{R}$, it suffices to define a quiver $\mathcal{Q}$ by taking $\mathcal{Q}_0$ as the set of objects in $\mathbf{R}$, while the arrows in $\mathcal{Q}_0$ will be those non-identity morphisms in $\mathbf{R}$ that do not admit any non-trivial factorization in $\mathbf{R}$. It is straightforward to check that $\mathbf{R}=\mathbf{Q}$.
\end{proof}
\end{proposition}

Proposition~\ref{ffdrc} may seem not completely satisfactory, because the degree functions do not necessarily match. The point is that any free direct Reedy category can be endowed with many different degree functions. In fact, Example~\ref{fdrc} describes just one of the possible degree functions on the free direct Reedy category associated to a quiver.

Given a small category $\mathbf{C}$ and an object $c\in\mathbf{C}$ we can consider the category $\mathbf{C}/c$ of arrows \emph{over} $c$. Namely, the objects of $\mathbf{C}/c$ are morphisms of $\mathbf{C}$ whose target is $c$. A morphism between two objects $s\to c$ and $t\to c$ in $\mathbf{C}/c$ is simply a commutative diagram of the form
\[ \xymatrix{  s\ar@{->}[dr]\ar@{->}[rr] & & t\ar@{->}[dl] \\
& c & .}  \]
Dually, one defines the category $c/\mathbf{C}$ of arrows \emph{under} $c$, where the objects are the morphisms in $\mathbf{C}$ whose source is $c$.

\begin{definition}\label{latching}
Suppose $\mathbf{R}$ is a free direct Reedy category and let $\mathbf{C}$ be a cocomplete category (i.e. a category $\mathbf{C}$ with all small colimits). Given a functor $F\colon \mathbf{R}\to\mathbf{C}$ and an object $r\in\mathbf{R}$, we define the \emph{latching object} at $r$ to be
\[ \Lt_r(F) = \colim\limits_{\mathbf{^{\downarrow}r}}F \]
where the colimit is taken over the full subcategory $\mathbf{^{\downarrow}r}$ of $\mathbf{R}/r$ containing all the objects except for the identity morphism $r\xrightarrow{\id_r} r$.
Dually, suppose $\mathbf{R}$ is a free inverse Reedy category and let $\mathbf{C}$ be a complete category (i.e. a category $\mathbf{C}$ with all small limits). Given a functor $F\colon \mathbf{R}\to\mathbf{C}$ and an object $r\in\mathbf{R}$, we define the \emph{matching object} at $r$ to be
\[ \Mt_r(F) = \lim\limits_{\mathbf{^{\uparrow}r}}F \]
where the limit is taken over the full subcategory $\mathbf{^{\uparrow}r}$ of $r/\mathbf{R}$ containing all the objects except for the identity morphism $r\xrightarrow{\id_r} r$. 
\end{definition}

We can now recall a powerful result that allows to ``lift" model structures to category of functors. Theorem~\ref{Reedy} holds in a much more general setting. For the strongest version we refer to \cite{Hir03}. We decided to present it in a weaker form in order to make the relation with our context as clear as possible.

\begin{theorem}\cite[Theorem 15.3.4]{Hir03}\label{Reedy}
Let $\mathbf{R}$ be a free direct Reedy category and let $\mathbf{M}$ be a model category. Then the category of functors $\mathbf{M}^{\mathbf{R}}$ admits a model structure where a natural transformation $F\to G$ is:
\begin{description}
\item[1] a Reedy weak equivalence if it is a pointwise weak equivalence in $\mathbf{M}$, i.e. if for every $r\in\mathbf{R}$ the morphism $F_r\to G_r$ is a weak equivalence in $\mathbf{M}$, 
\item[2] a Reedy-projective fibration if it is a pointwise fibration in $\mathbf{M}$, i.e. if for every $r\in\mathbf{R}$ the morphism $F_r\to G_r$ is a fibration in $\mathbf{M}$, 
\item[3] a Reedy-projective cofibration if for every $r\in\mathbf{R}$ the natural morphism induced by the pushout
\[ \Lt_r(G) \amalg_{\Lt_r(F)} F_r \longrightarrow G_r \]
is a cofibration in $\mathbf{M}$.
\end{description}
Dually, let $\mathbf{R}$ be a free inverse Reedy category and let $\mathbf{M}$ be a model category. Then the category of functors $\mathbf{M}^{\mathbf{R}}$ admits a model structure where a natural transformation $F\to G$ is:
\begin{description}
\item[1] a Reedy weak equivalence if it is a pointwise weak equivalence in $\mathbf{M}$, i.e. if for every $r\in\mathbf{R}$ the morphism $F_r\to G_r$ is a weak equivalence in $\mathbf{M}$, 
\item[2] a Reedy-injective cofibration if it is a pointwise cofibration in $\mathbf{M}$, i.e. if for every $r\in\mathbf{R}$ the morphism $F_r\to G_r$ is a cofibration in $\mathbf{M}$, 
\item[3] a Reedy-injective fibration if for every $r\in\mathbf{R}$ the natural morphism induced by the pullback
\[ F_r \longrightarrow \Mt_r(F) \times_{\Mt_r(G)} G_r \]
is a fibration in $\mathbf{M}$.
\end{description}
\end{theorem}

The model structures described in Theorem~\ref{Reedy} are usually called the Reedy-projective and the Reedy-injective model structures. We will often write $RP$-fibration (respectively $RP$-cofibration) instead of Reedy-projective fibration (respectively Reedy-projective cofibration). Similarly, we will refer to Reedy-injective fibrations (respectively Reedy-injective cofibrations) writing $RI$-fibrations (respectively $RI$-cofibrations).

\begin{remark}
The model structure described in Theorem~\ref{Reedy} is cofibrantly generated. In fact we restricted to the case where the Reedy category $\mathbf{R}$ is direct (respectively inverse), so that the Reedy model structure and the \emph{projective} (respectively \emph{injective}) \emph{model structure} on functors coincide, see \cite{Hir03}. Anyway, this is not going to play a crucial role in our contest.
\end{remark}

\begin{remark}
It is important to notice that the Reedy model structures described in Theorem~\ref{Reedy} do not depend on the degree function of the finite free direct (or inverse) Reedy category $\mathbf{R}$. That is, given two different degree functions $d$ and $d'$ on $\mathbf{R}$ the classes of fibrations, cofibrations and weak equivalences furnished by Theorem~\ref{Reedy} do not change.
\end{remark}

Our next goal is to describe the Reedy-projective cofibrant objects in the category $\mathbf{M}^{\mathbf{R}}$ with respect to the Reedy-projective model structure and, dually, the Reedy-injective fibrant objects in the category $\mathbf{M}^{\mathbf{R}}$ with respect to the Reedy-injective model structure.

\begin{lemma}\label{rc}
Given a model category $\mathbf{M}$ and a finite free direct Reedy category $\mathbf{R}$, a diagram $F\in\mathbf{M}^{\mathbf{R}}$ is $RP$-cofibrant if and only if for every $r\in\mathbf{R}$ the natural morphism
\[ \Lt_r(F)\to F_r \]
is a cofibration in $\mathbf{M}$.
Dually, given a model category $\mathbf{M}$ and a finite free inverse Reedy category $\mathbf{R}$, a diagram $F\in\mathbf{M}^{\mathbf{R}}$ is $RI$-fibrant if and only if for every $r\in\mathbf{R}$ the natural morphism
\[ F_r\to \Mt_r(F) \]
is a fibration in $\mathbf{M}$.
\begin{proof}
The statement immediately follows from Theorem~\ref{Reedy}.
\end{proof}
\end{lemma}

Suppose $\mathbf{R}$ is a finite free direct Reedy category. For every $r\in\mathbf{R}$ consider the category $\mathbf{^{\downarrow}r}$, see Definition~\ref{latching}. We now characterize the Reedy cofibrant objects in terms of the \emph{irreducible objects} of $\mathbf{^{\downarrow}r}$, that is the subset $\Irr(\mathbf{^{\downarrow}r})\subseteq\Ob(\mathbf{^{\downarrow}r})$ of arrows $s\to r$ which are irreducible morphisms in $\mathbf{R}$. Proposition~\ref{ffdrc} allows to give a ``quiver-theoretic" description of irreducible objects in $\mathbf{^{\downarrow}r}$. Suppose $\mathbf{R} = \mathbf{Q}$ for some finite acyclic quiver $\mathcal{Q}$ (up to the degree function). Then for every $r\in\mathbf{R}$ we have $\Irr(\mathbf{^{\downarrow}r}) = \{\alpha\in\mathcal{Q}_1 \,\,\vert\,\, \tau(\alpha)=r\}$. Roughly speaking, $\Irr(\mathbf{^{\downarrow}r})$ is the set of \emph{incoming} arrows at $r$. Dually, one can consider the subset $\Irr(\mathbf{^{\uparrow}r})\subseteq\Ob(\mathbf{^{\uparrow}r})$, which essentially represents the subset of \emph{outgoing} arrows. This observations together with Theorem~\ref{Rcof} will allow us to endow the category of quiver representations with explicit model structures, see Theorem~\ref{modelrep}.

\begin{theorem}\label{Rcof}
Given a model category $\mathbf{M}$ and a finite free direct Reedy category $\mathbf{R}$, a diagram $F\in\mathbf{M}^{\mathbf{R}}$ is $RP$-cofibrant if and only if for every $r\in\mathbf{R}$ the natural morphism
\[ \bigoplus\limits_{(s\to r)\in \Irr(\mathbf{^{\downarrow}r})}F_s \to F_r \]
is a cofibration in $\mathbf{M}$, where the direct sum is taken over all the objects in $\mathbf{^{\downarrow}r}$ that are irreducible morphisms in $\mathbf{R}$.
Dually, given a model category $\mathbf{M}$ and a finite free inverse Reedy category $\mathbf{R}$, a diagram $F\in\mathbf{M}^{\mathbf{R}}$ is $RI$-fibrant if and only if for every $r\in\mathbf{R}$ the natural morphism
\[ F_r \to \prod\limits_{(s\to r)\in \Irr(\mathbf{^{\uparrow}r})}F_s \]
is a fibration in $\mathbf{M}$, where the product is taken over all the objects in $\mathbf{^{\uparrow}r}$ that are irreducible morphisms in $\mathbf{R}$.
\begin{proof}
By Lemma~\ref{rc} it suffices to show that for every object $F\in\mathbf{M}^{\mathbf{R}}$ and for every $r\in\mathbf{R}$ there is an isomorphism $\Lt_r(F) \cong \bigoplus\limits_{(s\to r)\in \Irr(\mathbf{^{\downarrow}r})}F_s$. Let us fix $r\in\mathbf{R}$. An object in the category $\mathbf{^{\downarrow}r}$ is a path of the form $(i_n\to \cdots \to i_0\to r)$. Clearly, in $\mathbf{^{\downarrow}r}$ there exists a morphism from $(i_n\to \cdots \to i_0\to r)$ to $(i_0\to r)$. Since $\mathbf{R}$ is a free category, for every $(s\to r)\in \Irr(\mathbf{^{\downarrow}r})$ we have a ``connected component" $\mathbf{C}_s\subseteq \mathbf{^{\downarrow}r}$ defined as the full subcategory of $\mathbf{^{\downarrow}r}$ whose objects are the paths with an arrow to $(s\to r)$. Moreover, $(s\to r)$ is the final object of $\mathbf{C}_s$ for every ${(s\to r)\in \Irr(\mathbf{^{\downarrow}r})}$. Hence we have the following:
\[ \Lt_r F =\colim_{\mathbf{^{\downarrow}r}} F \cong \colim_{(s\to r)\in \Irr(\mathbf{^{\downarrow}r})}\left( \colim_{(i\to r)\in C_s}F_i \right) = \colim_{(s\to r)\in \Irr(\mathbf{^{\downarrow}r})}F_s = \bigoplus_{(s\to r)\in \Irr(\mathbf{^{\downarrow}r})}F_s. \]
The second part of the statement is dual.
\end{proof}
\end{theorem}

\begin{remark}
To understand the usefulness of Theorem~\ref{Reedy} and Theorem~\ref{Rcof} in the contest of quiver representations it is sufficient to observe the following: representations of a finite acyclic quiver $\mathcal{Q}$ over a quasi-Frobenius ring $R$ are precisely functors $\mathbf{Q} \to \Mod(R)$ where $\mathbf{Q}$ is the free direct Reedy category of Example~\ref{fdrc}, i.e. $\Rep(\mathcal{Q}, R) \cong \Mod(R)^{\mathbf{Q}}$. In fact, a functor $M\colon\mathbf{Q} \to \Mod(R)$ is the data of a module $M_j\in\Mod(R)$ for every $j\in\mathcal{Q}_0$, together with a $R$-linear map $M_{\sigma(\alpha)}\to M_{\tau(\alpha)}$ for every $\alpha\in\mathcal{Q}_1$.
\end{remark}

\begin{theorem}[Model structures on quiver representations]\label{modelrep}
Let $\mathcal{Q}$ be a finite acyclic quiver, and let $R$ be a quasi-Frobenius ring. Then the category $\Rep(\mathcal{Q}, R)$ admits two model structures.

In the \emph{Reedy projective} model structure a morphism between $R$-representations $M\to N$ is:
\begin{description}
\item[1] a weak equivalence (called Reedy stable equivalence) if it is a pointwise stable equivalence in $\Mod(R)$, i.e. if for every vertex $j\in \mathcal{Q}_0$ the morphism $M_j\to N_j$ is a stable equivalence of $R$-modules, 
\item[2] a fibration (called $RP$-fibration) if it is a pointwise surjection in $\Mod(R)$, i.e. if for every vertex $j\in \mathcal{Q}_0$ the morphism $M_j\to N_j$ is surjective, 
\item[3] a cofibration (called $RP$-cofibration) if for every vertex $j\in \mathcal{Q}_0$ the natural morphism induced by the coproduct
\[ \left( \bigoplus_{\stackrel{\alpha\in \mathcal{Q}_1}{\tau(\alpha) = j}}N_{\sigma(\alpha)} \right) \amalg_{\left( \bigoplus\limits_{\stackrel{\alpha\in \mathcal{Q}_1}{\tau(\alpha) = j}} M_{\sigma(\alpha)}\right)} M_j \longrightarrow N_j \]
is injective.
\end{description}
Moreover, every $R$-representation is $RP$-fibrant while an $R$-representation $M$ is $RP$-cofibrant if and only if for every vertex $j\in \mathcal{Q}_0$ the natural morphism of $R$-modules
\[ \bigoplus\limits_{\stackrel{\alpha\in \mathcal{Q}_1}{\tau(\alpha) = j}}M_{\sigma(\alpha)} \to M_j \]
is injective.

In the \emph{Reedy injective} model structure a morphism between $R$-representations $M\to N$ is:
\begin{description}
\item[1] a weak equivalence (called, again, Reedy stable equivalence) if it is a pointwise stable equivalence in $\Mod(R)$, i.e. if for every vertex $j\in \mathcal{Q}_0$ the morphism $M_j\to N_j$ is a stable equivalence of $R$-modules, 
\item[2] a cofibration (called $RI$-cofibration) if it is a pointwise injection in $\Mod(R)$, i.e. if for every vertex $j\in \mathcal{Q}_0$ the morphism $M_j\to N_j$ is injective, 
\item[3] a fibration (called $RI$-fibration) if for every vertex $j\in \mathcal{Q}_0$ the natural morphism induced by the product
\[ M_j \longrightarrow \left( \prod_{\stackrel{\alpha\in \mathcal{Q}_1}{\sigma(\alpha) = j}}M_{\tau(\alpha)} \right) \times_{\left( \prod\limits_{\stackrel{\alpha\in \mathcal{Q}_1}{\sigma(\alpha) = j}} N_{\tau(\alpha)}\right)} N_j \]
is surjective.
\end{description}
Moreover, every $R$-representation is $RI$-cofibrant, while an $R$-representation $M$ is $RI$-fibrant if and only if for every vertex $j\in \mathcal{Q}_0$ the natural morphism of $R$-modules
\[ N_j\to \prod\limits_{\stackrel{\alpha\in \mathcal{Q}_1}{\sigma(\alpha) = j}}M_{\tau(\alpha)} \]
is surjective.
\begin{proof}
Since there exists an isomorphism of categories $\Rep(\mathcal{Q}, R)\cong \Mod(R)^{\mathbf{Q}}$, the statement immediately follows from Theorem~\ref{Reedy} and Theorem~\ref{Rcof}.
\end{proof}
\end{theorem}

Morphisms of $\Rep(\mathcal{Q},R)$ that are both Reedy stable equivalences and $RP$-fibrations (respectively $RP$-cofibrations) will be called $RP$-stable fibrations (respectively $RP$-stable cofibrations).
Morphisms of $\Rep(\mathcal{Q},R)$ that are both Reedy stable equivalences and $RI$-fibrations (respectively $RI$-cofibrations) will be called $RI$-stable fibrations (respectively $RI$-stable cofibrations).

\begin{remark}\label{generalization}
In the introduction we mentioned that our arguments work in a more general setting. More precisely one may assume $R$ to be a possibly non-commutative Gorenstein ring, see Definition~\ref{iwanagagorenstein}. In particular, Theorem~\ref{modelrep} would give a different characterization of $RP$-cofibrant $R$-representations, namely the natural morphism
\[ \bigoplus\limits_{\stackrel{\alpha\in \mathcal{Q}_1}{\tau(\alpha) = j}}M_{\sigma(\alpha)} \to M_j \]
is required to be a cofibration in $\Mod(R)$, i.e. injective with Gorenstein-projective cokernel. Clearly, if $R$ is quasi-Frobenius the second condition is automatically satisfied. We point out that even in the general case of a non-commutative Gorenstein ring $R$, given an $RP$-cofibrant $R$-representation $M\in\Rep(\mathcal{Q}, R)$ then $M_j$ is Gorenstein-projective in $\Mod(R)$ for every $j\in\mathcal{Q}_0$. In fact one can easily show that any $RP$-cofibration is in particular a vertexwise cofibration in $\Mod(R)$. As a consequence we have that the ``monic representations satisfying condition (G)" introduced in \cite{LZ13} are nothing but $RP$-cofibrant representations of Theorem~\ref{modelrep}.
\end{remark}

\begin{remark}
In \cite{LZ13} the approach via monic representations works only in the ``finite" case, i.e. the authors restrict their study to finitely generated modules. On the other hand, a more recent paper by Luo and Zhang extends their results to quivers with relations, see~\cite{LZ17}. Since model categories are not used very often in quiver representation theory, we decided to reduce technical issues at minimum. This motivates our choice to avoid the study of quivers with relations. Nevertheless, our approach via model categories works for possibly infinitely generated modules. However, the strategy of the present paper should give an explicit characterization of Gorenstein-projective modules even in the general case of infinitely generated Gorenstein-projective modules over the (opposite) path algebra of quivers with relations.

We point out that recently, Hu, Luo, Xiong and Zhou studied Gorenstein-projective modules in the non-finite case via filtration categories, see~\cite{HLXZ}.
Moreover, two papers by Enochs, Estrada, Garc\'ia Rozas,~\cite{EEG09}, and by Eshraghi, Hafezi, Salarian,~\cite{EHS13}, completely solved the problem providing a description of Gorenstein-projective modules in a very general setting, so that the main innovation of the present paper relies in the homotopic approach which is far from the ones explored up to now.
\end{remark}

\section{Stable equivalences and Reedy stable equivalences}\label{searse}

In order to better understand the notion of Reedy stable equivalences of Theorem~\ref{modelrep} we now look for a description as simple as possible for the class of stable equivalences in $\Mod(R)$, when $R$ is a quasi-Frobenius ring. This is already known; for the reader convenience we provide elementary proofs, see Corollary~\ref{stableeq}.

\begin{proposition}\label{stablecoffib}
Let $R$ be a quasi-Frobenius ring. A morphism $f\colon M\to N$ in $\Mod(R)$ is a stable fibration if and only if there exists an injective/projective $R$-module $P$ such that $f$ is the composition $M\xrightarrow{\cong} N\oplus P \to N$, where $N\oplus P\to N$ is the natural projection and $M\xrightarrow{\cong} N\oplus P$ is an isomorphism of $R$-modules.

Dually, a morphism $f\colon M\to N$ in $\Mod(R)$ is a stable cofibration if and only if there exists an injective/projective $R$-module $J$ such that $f$ is the composition $M\hookrightarrow M\oplus J \xrightarrow{\cong} N$, where $M\hookrightarrow M\oplus J$ is the natural inclusion and $M\oplus J\xrightarrow{\cong} N$ is an isomorphism of $R$-modules.
\begin{proof}
Suppose $f\colon M\to N$ in $\Mod(R)$ is a stable fibration. By Theorem~\ref{Mark} it is surjective and hence $f$ fits into a short exact sequence
\[ 0\to\ker(f)\xrightarrow{\iota} M\xrightarrow{f} N\to 0 \]
in $\Mod(R)$. Now, $N$ is cofibrant and then the diagram of solid arrows
\[ \xymatrix{ & M\ar@{->>}[d]_{\sim}^{f} \\
N\ar@{->}[r]_{\id_N} \ar@{.>}[ur]^{h} & N } \]
admits a lifting $h\colon N\to M$. This means that the above short exact sequence splits, i.e. $f$ is the composition $M\xrightarrow{\cong}\ker(f)\oplus N\to N$. It remains to show that $\ker(f)$ is injective/projective. Notice that the splitting gives a morphism $M\xrightarrow{r} \ker(f)$ such that $r\circ\iota = \id_{\ker(f)}$. Also, being $f$ a stable equivalence, $\ker(f)\xrightarrow{\iota} M$ is zero in $\underline{\Mod}(R)$, that is $\iota$ factors through a projective $R$-module in $\Mod(R)$. Now, the relation $r\circ\iota = \id_{\ker(f)}$ implies that $\id_{\ker(f)}$ factors through a projective object in $\Mod(R)$, and then $\ker(f)$ is a projective $R$-module as required. The proof of the second part of the statement is dual.
\end{proof}
\end{proposition}

\begin{corollary}\label{stableeq}
Let $R$ be a quasi-Frobenius ring. Every stable equivalence $M\xrightarrow{f} N$ in $\Mod(R)$ is of the form $M\to M\oplus P\cong N\oplus Q \to N$ for some injective/projective $R$-modules $P$ and $Q$, where $M\to M\oplus P$ is the natural inclusion and $N\oplus Q \to N$ is the natural projection.
\begin{proof}
Consider a stable equivalence $M\xrightarrow{f} N$ in $\Mod(R)$. Now, take a factorization as a cofibration followed by a stable fibration $f\colon M\xrightarrow{\alpha} A\xrightarrow{\beta} N$. By the two-out-of-three axiom $\alpha$ is a stable cofibration, then Proposition~\ref{stablecoffib} concludes the proof.
\end{proof}
\end{corollary}

The next step is to investigate $RP$-stable fibrations and $RI$-stable cofibrations in the category $\Rep(\mathcal{Q}, R)$.

\begin{proposition}\label{Rstablefib}
Let $\mathcal{Q}$ be a finite acyclic quiver, and let $R$ be a quasi-Frobenius ring. A morphism $f\colon M\to N$ in $\Rep(\mathcal{Q}, R)$ is a $RP$-stable fibration with respect to the model structure of Theorem~\ref{modelrep} if and only if for every vertex $j\in\mathcal{Q}_0$ the morphism $f_j\colon M_j\to N_j$ is a stable fibration in $\Mod(R)$.

Dually, a morphism $f\colon M\to N$ in $\Rep(\mathcal{Q}, R)$ is a $RI$-stable cofibration with respect to the model structure of Theorem~\ref{modelrep} if and only if for every vertex $j\in\mathcal{Q}_0$ the morphism $f_j\colon M_j\to N_j$ is a stable cofibration in $\Mod(R)$.
\begin{proof}
The statement is a corollary of Theorem~\ref{modelrep} and Proposition~\ref{stablecoffib}.
\end{proof}
\end{proposition}

Proposition~\ref{Rstablefib} may seem not completely satisfactory because the factorization is only degreewise in the category $\Mod(R)$ and \emph{not} global in the category $\Rep(\mathcal{Q},R)$. Namely one could naturally define $RP$-elementary stable fibrations in $\Rep(\mathcal{Q},R)$ as the compositions $M\xrightarrow{\cong} N\oplus P \xrightarrow{\pi_N} N$ for some projective $R$-representation $P\in\Rep(Q,R)$, and may expect that $RP$-stable fibrations coincides with $RP$-elementary stable fibrations generalizing Proposition~\ref{stablecoffib}. (Un)fortunately this turns out to be false, and we will see in Section~\ref{gorenstein} how this asymmetry naturally gives rise to the notion of Gorenstein-projective modules. In fact, every diagram of solid arrows of the following shape
\[ \xymatrix{ & M\ar@{->}[d]^{\cong} \\
 & N\oplus P\ar@{->}[d]^{\pi_N} \\
C\ar@{->}[r]\ar@{.>}[uur]^{h} & N } \]
admits a lifting $h$ even in the category $\Rep(\mathcal{Q},R)$. This means that if $RP$-elementary stable fibrations were precisely the $RP$-stable fibrations, then every $R$-representation would be $RP$-cofibrant but this is clearly false by Theorem~\ref{modelrep}, unless the quiver has no arrows.

Our next result characterizes projective and injective $R$-representations in terms of the Reedy model structures of Theorem~\ref{modelrep}.

\begin{lemma}\label{projrep}
Let $\mathcal{Q}$ be a finite acyclic quiver, and let $R$ be a quasi-Frobenius ring. Then the following are equivalent for an $R$-representation $M$.
\begin{description}
\item[1] $M$ is a projective object in $\Rep(\mathcal{Q},R)$.
\item[2] $M$ is $RP$-cofibrant and $M_j$ is a projective $R$ module for every $j\in\mathcal{Q}_0$.
\end{description}
Dually, the following are equivalent.
\begin{description}
\item[1] $M$ is an injective object in $\Rep(\mathcal{Q},R)$.
\item[2] $M$ is $RI$-fibrant and $M_j$ is an injective $R$ module for every $j\in\mathcal{Q}_0$.
\end{description}
\begin{proof}
Suppose $M$ is a projective $R$-representation. Then it is clearly cofibrant since it satisfies the left lifting property with respect to every surjective morphisms and then, in particular, with respect to every $RP$-stable fibration. It remains to show that $M_j$ is projective in $\Mod(R)$ for every $j\in\mathcal{Q}_0$. To this aim, we fix $i\in\mathcal{Q}_0$ and given a diagram of solid arrows of shape
\[ \xymatrix{  & A\ar@{->>}[d]^{p} \\
M_i\ar@{->}[r] \ar@{.>}[ur]^{h} & B } \]
in the category $\Mod(R)$ we construct a (dotted) lifting $h\colon M_i\to A$. We begin by defining two new $R$-representations $A(i), B(i)\in\Rep(\mathcal{Q},R)$ through a general procedure. For any $C\in\Mod(R)$, define
\[ C(i)_j = \prod_{\text{paths} \;  j\to i} C \]
and for any arrow $\alpha\colon j\to j'$ let the corresponding morphism $C(i)_{\alpha} \colon   C(i)_j \to  C(i)_{j'}$ be the canonical projection, defined mapping (via the identity) the component $C$ corresponding to a path $j\xrightarrow{\alpha} j' \to i$ to the component $C$ corresponding to the path $j'\to i$, and being $0$ on any component corresponding to a path $j\to i$ which does not factorize through $\alpha$. In particular, $A(i)_i=A$. The above procedure defines an exact functor
 \[ (i)\colon \Mod(R) \to \Rep(\mathcal{Q},R), \qquad C\mapsto C(i) \]
so that the epimorphism of $R$-modules $p: A \to B$ induces an obvious epimorphism of $R$-representations $p(i): A(i) \to B(i)$. Moreover, this functor is right adjoint to the projection functor $M\mapsto M_i$. Hence the map $M_i \to B$ naturally induces a morphism $M\to B(i)$ which by hypothesis can be lifted to a morphism $\tilde{h}\colon M \to A(i)$ such that the following diagram
\[ \xymatrix{  & A(i)\ar@{->}[d]^{p(i)} \\
M\ar@{->}[r] \ar@{.>}[ur]^{\tilde{h}} & B(i) } \]
commutes in the category $\Rep(\mathcal{Q},R)$. Again by adjunction, it is now sufficient to define $h = \tilde{h}_i\colon M_i\to A$ to complete the proof.

Viceversa, given a diagram of solid arrows of shape
\[ \xymatrix{  & X\ar@{->>}[d]^{p} \\
M\ar@{->}[r]_{g} \ar@{.>}[ur]^{f} & Y } \]
in the category $\Rep(\mathcal{Q},R)$, we can recursively construct the (dotted) lifting $f$ using the injectivity of the natural morphisms
\[ \bar{m}_i\colon  \bigoplus\limits_{\stackrel{\alpha\in\mathcal{Q}_1}{\tau(\alpha)=i}}M_{\sigma(\alpha)}\to M_i \]
for every $i\in\mathcal{Q}_0$, and the degreewise projectivity of $M$. Let us be more precise. First, if $i\in\mathcal{Q}_0$ is a source vertex (i.e. it does not exist any arrow whose target is $i$), then we can define $f_i\colon M_i\to X_i$ to be any map of $R$-modules satisfying $g_i=p_if_i$, which exists by the projectivity of $M_i$. Now, suppose $i\in \mathcal{Q}_0$ is \emph{any} vertex and assume we have already defined $f_j\colon M_j\to X_j$ for every $j\in\mathcal{Q}_0$ such that there exists a path $j\to j_1\to \cdots \to j_k\to i$. The universal property of the direct sum gives a map
\[ \bar{f} \colon \bigoplus\limits_{\stackrel{\alpha\in\mathcal{Q}_1}{\tau(\alpha)=i}}M_{\sigma(\alpha)}\xrightarrow{\oplus f_{\sigma(\alpha)}} X_i \]
satisfying the following commutativity
\begin{equation}\label{eq.lemma2.4}
\xymatrix{	M_j \ar@{->}[rr]^{f_j} \ar@{->}[d]_{\varphi_i} & & X_j \ar@{->}[d]^{x_{ji}} \\
\bigoplus\limits_{\stackrel{\alpha\in\mathcal{Q}_1}{\tau(\alpha)=i}}M_{\sigma(\alpha)} \ar@{->}[rr]^{\bar{f}} \ar@{->}[d]^{\bar{m}_i} & & X_i \ar@{->}[d]^{p_i} \\
M_i \ar@{->}[rr]^{g_i} & & Y_i	}
\end{equation}
for every arrow $j\to i$, where $\varphi_i$ is the natural inclusion. Now observe that the map $\bar{m}_i$ is injective being $M$ an $RP$-cofibrant $R$-representation by hypothesis, and $p_i$ is surjective. Moreover, the $R$-module $\bigoplus\limits_{\stackrel{\alpha\in\mathcal{Q}_1}{\tau(\alpha)=i}}M_{\sigma(\alpha)}$ is projective (hence injective being $R$ a quasi-Frobenius ring) so that the short exact sequence
\[ 0\to \bigoplus\limits_{\stackrel{\alpha\in\mathcal{Q}_1}{\tau(\alpha)=i}}M_{\sigma(\alpha)} \xrightarrow{\bar{m}_i} M_i \to \coker(\bar{m}_i) \to 0 \]
splits, and $\coker(\bar{m}_i)$ is a projective $R$-module being isomorphic to a direct summand of $M_i$. It follows by Proposition~\ref{stablecoffib} that the map $\bar{m}_i$ is a stable cofibration in $\Mod(R)$, so that the bottom square in~\eqref{eq.lemma2.4} admits a lifting $f_i\colon M_i\to X_i$ as required. Hence $M$ is projective in $\Rep(\mathcal{Q},R)$.
The second part of the statement is dual.
\end{proof}
\end{lemma}

%The author is grateful to an anonymous referee for useful comments and for requiring more details making the proof of Lemma~\ref{projrep} more readable and precise.
We now give another description of $RP$-stable fibrations and $RI$-stable cofibrations.

\begin{proposition}\label{finiteRTF}
Let $\mathcal{Q}$ be a finite acyclic quiver, and let $R$ be a quasi-Frobenius ring. A morphism $f\colon M\to N$ in $\Rep(\mathcal{Q}, R)$ is a $RP$-stable fibration if and only if it is surjective and $\projdim_{\Rep(\mathcal{Q},R)}\left(\ker\{f\}\right)\leq1$.

Dually, a morphism $f\colon M\to N$ in $\Rep(\mathcal{Q}, R)$ is a $RI$-stable cofibration if and only if it is injective and $\injdim_{\Rep(\mathcal{Q},R)}\left(\coker\{f\}\right)\leq1$.
\begin{proof}
Let $f\colon M\to N$ in $\Rep(\mathcal{Q}, R)$ be a $RP$-stable fibration. In order to show that
\[ \projdim_{\Rep(\mathcal{Q},R)}\left(\ker\{f\}\right)\leq1 \]
we construct an explicit projective resolution of $\ker\{f\}$ in the category $\Rep(\mathcal{Q},R)$. To begin with, for every $j\in\mathcal{Q}_0$ consider the category $\mathbf{^{\downarrow}j}$ of paths ending at $j$, the notation is the same as in Definition~\ref{latching}. Recall that by definition the identity path $j\xrightarrow{\id} j$ does not belong to $\mathbf{^{\downarrow}j}$, so that it may be empty. We denote by $\overline{\mathbf{^{\downarrow}j} }$ the ``closure" of $\mathbf{^{\downarrow}j}$, namely $\overline{\mathbf{^{\downarrow}j} } = \mathbf{^{\downarrow}j} \cup \{j\xrightarrow{\id} j\}$. In order to keep the exposition as clear as possible, we shall denote by $\sigma(\beta)\in\mathcal{Q}_0$ the starting vertex for a path $\beta\in\overline{\mathbf{^{\downarrow}j}}$; of course $\tau(\beta)=j$ for every $\beta\in\overline{\mathbf{^{\downarrow}j}}$.
We now define an $R$-representation $T\in\Rep(\mathcal{Q},R)$ as follows. For every vertex $j\in\mathcal{Q}_0$ we set $T_j=\bigoplus\limits_{\beta\in\overline{\mathbf{^{\downarrow}j}}}\ker\{f_{\sigma(\beta)}\}$. Clearly, $\overline{\mathbf{^{\downarrow}\sigma(\alpha)}}$ is a subcategory of $\overline{\mathbf{^{\downarrow}\tau(\alpha)}}$ for every arrow $\alpha\in\mathcal{Q}_1$, then there exists a natural inclusion
\[ \overline{\iota}\colon\bigoplus\limits_{\beta\in\overline{\mathbf{^{\downarrow}\sigma(\alpha)}}}\ker\{f_{\sigma(\beta)}\}\to\bigoplus\limits_{\beta\in\overline{\mathbf{^{\downarrow}\tau(\alpha)}}}\ker\{f_{\sigma(\beta)}\} \]
so that we can define
\[ t_{\alpha}=\overline{m}_{\alpha} + \overline{\iota}\colon T_{\sigma(\alpha)}\to T_{\tau(\alpha)} \]
where $\overline{m}_{\alpha}\colon\ker\{f_{\sigma(\alpha)}\}\to\ker\{f_{\tau(\alpha)}\}$ is the restriction of $m_{\alpha}\colon M_{\sigma(\alpha)}\to M_{\tau(\alpha)}$ to the submodule $\ker\{f_{\sigma(\alpha)}\}\subseteq M_{\sigma(\alpha)}$. Since for every $j\in\mathcal{Q}_0$ the multiplicity of $\ker\{f_j\}$ in $T_j$ is precisely $1$, there exists an obvious projection $\pi\colon T\to \ker\{f\}$ in the category $\Rep(\mathcal{Q},R)$. Moreover, there is a short exact sequence in the category $\Rep(\mathcal{Q},R)$
\[ 0\to K\to T\xrightarrow{\pi} \ker\{f\}\to 0 \]
where for every vertex $j\in\mathcal{Q}_0$ we define $K_j=\bigoplus\limits_{\beta\in\mathbf{^{\downarrow}j}}\ker\{f_{\sigma(\beta)}\}$, and for every $\alpha\in\mathcal{Q}_1$ the morphism $k_{\alpha}\colon K_{\sigma(\alpha)} \to K_{\tau(\alpha)}$ is the natural inclusion induced by $\mathbf{^{\downarrow}\sigma(\alpha)}\subseteq\mathbf{^{\downarrow}\tau(\alpha)}$. It now suffices to show that $K$ and $T$ are projective objects in the category $\Rep(\mathcal{Q},R)$. By Proposition~\ref{Rstablefib}, for every vertex $j\in\mathcal{Q}_0$ the morphism $f_j\colon M_j\to N_j$ is an elementary stable equivalence in $\Mod(R)$, i.e. $f_j$ is the composition $f_j\colon M_j\xrightarrow{\cong} N_j\oplus P_j\to N_j$ for some projective $R$-modules $P_j\in\Mod(R)$. It follows that $\ker\{f_j\}\cong P_j$ in $\Mod(R)$ for every $j\in\mathcal{Q}_0$, so that $K_j$ and $T_j$ are projective $R$-modules. Then, the statement follows by Lemma~\ref{projrep} using the injectivity of the natural morphisms
\[ \bigoplus\limits_{\stackrel{\alpha\in\mathcal{Q}_1}{\tau(\alpha)=j}}K_{\sigma(\alpha)}\to K_j \qquad \bigoplus\limits_{\stackrel{\alpha\in\mathcal{Q}_1}{\tau(\alpha)=j}}T_{\sigma(\alpha)}\to T_j \]
for every $j\in\mathcal{Q}_0$.

For the converse, let $f\colon M\to N$ be a surjective morphism in $\Rep(\mathcal{Q},R)$ such that $\ker\{f\}$ satisfies $\projdim_{\Rep(\mathcal{Q},R)}(\ker\{f\})\leq 1$. Then there exists a short exact sequence
\[ 0\to K\to T\to \ker\{f\}\to 0 \]
with $K$ and $T$ projective objects in $\Rep(\mathcal{Q},R)$. Since for every $j\in\mathcal{Q}_0$ the $R$-module $K_j$ is (projective and) injective, the sequence
\[ 0\to K_j\to T_j\to \ker\{f_j\}\to 0 \]
splits in $\Mod(R)$. It follows that $\ker\{f_j\}$ is isomorphic to a direct summand of $T_j$, and then it is a projective $R$-module. Now, for every $j\in\mathcal{Q}_0$ we consider the short exact sequence
\[ 0\to \ker\{f_j\}\to M_j\xrightarrow{f_j} N_j \to 0. \]
Since $R$ is quasi-Frobenius, $\ker\{f_j\}$ is an injective $R$-module so that the sequence above splits in $\Mod(R)$. Hence by Proposition~\ref{stablecoffib} $f_j$ is a stable fibration as required.
The second part of the statement is dual.
\end{proof}
\end{proposition}

The last result of this section is an easy consequence of Proposition~\ref{finiteRTF}.

\begin{corollary}
Let $\mathcal{Q}$ be a finite acyclic quiver, and let $R$ be a quasi-Frobenius ring. For every $R$-representation $M\in\Rep(\mathcal{Q}, R)$ there exists an exact sequence
\[ 0\to S\to T\to G_PM\to M\to 0 \]
where $S$ and $T$ are projective objects in $\Rep(\mathcal{Q},R)$ and $G_PM$ is $RP$-cofibrant.

Dually, there exists an exact sequence
\[ 0\to M\to G_IM\to I\to J\to 0 \]
where $I$ and $J$ are injective objects in $\Rep(\mathcal{Q},R)$ and $G_IM$ is a $RI$-fibrant $R$-representation.
\begin{proof}
Take a factorization of the morphism $0\to M$ as a $RP$-cofibration followed by a $RP$-stable fibration: $0\to G_PM\xrightarrow{\epsilon} M$. By Proposition~\ref{finiteRTF} the kernel $K=\ker\{\epsilon\}$ has projective dimension at most $1$, so that in $\Rep(\mathcal{Q},R)$ there exists an exact sequence
\[ 0\to S\to T\xrightarrow{f} K\to 0 \]
with $S$ and $T$ projective. Now
\[ 0\to S\to T\xrightarrow{\iota\circ f} G_PM\xrightarrow{\epsilon} M\to 0 \]
is the required exact sequence, where $\iota\colon K\to G_PM$ is the natural inclusion. The second part of the statement is dual.
\end{proof}
\end{corollary}

\section{Relation with Gorenstein Homological Algebra}\label{gorenstein}

Throughout this section we will consider a finite acyclic quiver $\mathcal{Q}$, and the path algebra $\Lambda$ of $\mathcal{Q}$ over a quasi-Frobenius ring $R$.
Recall that $\Lambda$ is generated as an $R$-module by all paths in $\mathcal{Q}$ of length greater than or equal to zero, so that $\Lambda$ includes the so-called \emph{lazy paths}, one for each vertex of the quiver. The multiplication in $\Lambda$ is given by composition of paths, and if two paths cannot be concatenated, then by definition their product in $\Lambda$ is $0$. Notice that this defines an associative algebra over $R$. This algebra has a unit element since the quivers we are interested in are assumed to have only finitely many vertices. Historically, composition of paths is written from the left to the right. To avoid confusion, we shall write $\Lambda^{op}$ instead of $\Lambda$ when dealing with the opposite product convention. Following this notation, modules over $\Lambda^{op}$ are naturally identified with representations of $\mathcal{Q}$.
Moreover, under our assumptions on the quiver, if $R=\mathbb{K}$ is a field then $\Lambda$ is a finite-dimensional hereditary algebra over $R$, see \cite{Rin11}. In Lemma~\ref{1gor} we prove that $\Lambda$ is in fact $1$-Gorenstein for every quasi-Frobenius ring $R$.

Many interesting examples of such algebras arise when $R$ is a finite dimensional self-injective $\mathbb{C}$-algebra such as $R=\frac{\mathbb{C}[t]}{(t^n)}$, in this case we have $\Lambda = \mathbb{C}\mathcal{Q}\otimes_{\mathbb{C}} R$ where $\mathbb{C}\mathcal{Q}$ is the usual path algebra of $\mathcal{Q}$ over $\mathbb{C}$.

We begin by recalling the notions of $n$-\emph{Gorenstein} rings, \emph{Gorenstein-projective} modules and \emph{Gorenstein-injective} modules.

\begin{definition}\label{iwanagagorenstein}
Given $n\in\mathbb{N}$, a Noetherian ring $G$ is called $n$-\emph{Gorenstein} if
\[ \injdim_{\Mod(G)} G \leq n \mbox{ and } \injdim_{\Mod(G^{op})} G \leq n. \]
That is, $G$ has injective dimension at most $n$ both as a left and right module over itself. $G$ is called \emph{Gorenstein} if it is $n$-Gorenstein for some $n\in\mathbb{N}$.
\end{definition}

Gorenstein rings were introduced by Iwanaga in \cite{Iw79} and \cite{Iw80}, generalizing the standard notion of commutative Gorenstein rings. Examples of Gorenstein rings are quasi-Frobenius rings and group rings $K[G]$ for any commutative Gorenstein ring $K$ and any finite group $G$, see \cite{EN55}.

\begin{definition}
Let $R$ be a ring. An $R$-module $M\in\Mod(R)$ is called \emph{Gorenstein-projective} if there exists an exact sequence of projective modules
\[ \cdots\to P^{-1}\xrightarrow{d^{-1}} P^0\xrightarrow{d^0} P^1\to \cdots \]
that remains exact under the functor $\Hom(-,P)$ for every projective module $P\in\Mod(R)$, and such that $M\cong \ker\{d^0\}$.

Dually, $M\in\Mod(R)$ is called \emph{Gorenstein-injective} if there exists an exact sequence of injective modules
\[ \cdots\to J^{-1}\xrightarrow{d^{-1}} J^0\xrightarrow{d^0} J^1\to \cdots \]
that remains exact under the functor $\Hom(J,-)$ for every injective module $J\in\Mod(R)$, and such that $M\cong \ker\{d^0\}$.
\end{definition}

We shall denote by $\GProj(\Lambda)$ the full subcategory of $\Mod(\Lambda)$ whose objects are the Gorenstein-projective modules and, similarly, by $\GInj(\Lambda)$ the full subcategory of $\Mod(\Lambda)$ whose objects are the Gorenstein-injective modules. The main results concerning Gorenstein-projective modules are described by Enochs and Jenda, and can be found in \cite{EJ11}. Over a Gorenstein ring, finitely generated Gorenstein-projective modules are also called maximal Cohen-Macaulay modules. The relation between finitely generated Gorenstein-projective modules and quiver representations over a finite dimensional self-injective algebra is investigated by Luo and Zhang in \cite{LZ13}. We will extend one of their results to the whole class of Gorenstein-projective modules, see Corollary~\ref{gorcof}.

In \cite{Hov02}, Hovey introduced two model structures on the category $\Mod(G)$ for any Gorenstein ring $G$, and he studied the main properties of the associated homotopy category $\Ho(\Mod(G))$ with respect to this model structures. We will call these model structures \emph{Hovey-projective} (or simply $HP$) and \emph{Hovey-injective} (or $HI$) structures. Our next goal is to induce a model structure on the category $\Mod(\Lambda)$ through the equivalence of categories $\Rep(\mathcal{Q},R) \simeq \Mod(\Lambda)$, and to compare it with Hovey's model structures. To this aim we need to show that given a finite acyclic quiver $\mathcal{Q}$ and a quasi-Frobenius ring $R$, the path algebra $\Lambda = R\mathcal{Q}$ is a Gorenstein ring, see Lemma~\ref{1gor}. We begin by recalling Hovey's result.

\begin{theorem}[Hovey, \cite{Hov02}]\label{Hovey}
Suppose $G$ is a Gorenstein ring. Then there are two model structures on the category of $G$-modules where the class of trivial objects is the class of modules of finite projective dimension. In the Hovey-projective model structure
\begin{description}
\item[1] HP-fibrations coincide with surjective morphisms,
\item[2] trivial HP-fibrations are surjective morphisms whose kernels have finite projective dimension,
\item[3] HP-cofibrant objects are the Gorenstein-projective modules,
\item[4] modules that are both HP-cofibrant and trivial are projective.
\end{description}
In the Hovey-injective model structure
\begin{description}
\item[1] HI-cofibrations coincide with injective morphisms,
\item[2] trivial HI-cofibrations are injective morphisms whose cokernels have finite injective dimension,
\item[3] HI-fibrant objects are the Gorenstein-injective modules.
\item[4] modules that are both HI-fibrant and trivial are injective.
\end{description}
\end{theorem}

Now we state the existence of two model structures on $\Mod(\Lambda^{op})$, which we call again Reedy-projective and Reedy-injective model structures since they are induced by Theorem~\ref{modelrep}.

\begin{theorem}\label{modelmodule}
Let $\mathcal{Q}$ be a finite acyclic quiver, and let $R$ be a quasi-Frobenius ring. Consider the path algebra $\Lambda=R\mathcal{Q}$. Also, for every vertex $j\in \mathcal{Q}_0$ let us denote by $e_j\in\Lambda$ the correspondent lazy path. Then the category $\Mod(\Lambda^{op})$ of $\Lambda^{op}$-modules admits two model structures. In the \emph{Reedy-projective} model structure:
\begin{description}
\item[1] a morphism of $\Lambda^{op}$-modules $M\to N$ is a weak equivalence (called Reedy stable equivalence) if and only if for every vertex $j\in \mathcal{Q}_0$ the induced morphism $e_j(M)\to e_j(N)$ is a stable equivalence of $R$-modules,
\item[2] a morphism of $\Lambda^{op}$-modules is a $RP$-fibration if and only if it is surjective,
\item[3] a morphism of $\Lambda^{op}$-modules $M\to N$ is a $RP$-cofibration if for every vertex $j\in \mathcal{Q}_0$ the natural morphism induced by the pushout
\[ \left( \bigoplus_{\stackrel{\alpha\in \mathcal{Q}_1}{\tau(\alpha) = j}}e_{\sigma(\alpha)}(N) \right) \amalg_{\left( \bigoplus\limits_{\stackrel{\alpha\in \mathcal{Q}_1}{\tau(\alpha) = j}} e_{\sigma(\alpha)}(M)\right)} e_j(M) \longrightarrow e_j(N) \]
is injective.
\item[4] a $\Lambda^{op}$-module $M\in\Mod(\Lambda^{op})$ is $RP$-cofibrant if and only if the natural morphism
\[ \bigoplus_{\stackrel{\alpha\in \mathcal{Q}_1}{\tau(\alpha) = j}}e_{\sigma(\alpha)}(M) \to e_j(M) \]
is injective for every $j\in\mathcal{Q}_0$.
\end{description}
Dually, in the \emph{Reedy-injective} model structure:
\begin{description}
\item[1] a morphism of $\Lambda^{op}$-modules $M\to N$ is a weak equivalence (called Reedy stable equivalence) if and only if for every vertex $j\in \mathcal{Q}_0$ the induced morphism $e_j(M)\to e_j(N)$ is a stable equivalence of $R$-modules,
\item[2] a morphism of $\Lambda^{op}$-modules is a $RP$-cofibration if and only if it is injective,
\item[3] a morphism of $\Lambda^{op}$-modules $M\to N$ is a $RP$-fibration if for every vertex $j\in \mathcal{Q}_0$ the natural morphism induced by the pullback
\[ e_j(M)\longrightarrow \left( \prod_{\stackrel{\alpha\in \mathcal{Q}_1}{\tau(\alpha) = j}}e_{\tau(\alpha)}(M) \right) \times_{\left( \prod\limits_{\stackrel{\alpha\in \mathcal{Q}_1}{\sigma(\alpha) = j}} e_{\tau(\alpha)}(N)\right)} e_j(N) \]
is injective.
\item[4] a $\Lambda^{op}$-module $M\in\Mod(\Lambda^{op})$ is $RI$-fibrant if and only if the natural morphism
\[ e_j(M) \to \prod_{\stackrel{\alpha\in \mathcal{Q}_1}{\sigma(\alpha) = j}}e_{\tau(\alpha)}(M) \]
is surjective for every $j\in\mathcal{Q}_0$.
\end{description}
\begin{proof}
The statement is essentially the same as the one of Theorem~\ref{modelrep}. First of all recall that there exists an equivalence of categories $\Mod(\Lambda^{op})\simeq\Rep(\mathcal{Q},A)$, which assigns to every $\Lambda^{op}$-module $M$ the $R$-representation whose $R$-module over a vertex $j\in \mathcal{Q}_0$ is $e_j(M)$, while for every $\alpha\in \mathcal{Q}_1$ the $R$-linear morphism $e_{\sigma(\alpha)}(M)\to e_{\tau(\alpha)}(M)$ is given by the action of $\alpha\in\Lambda^{op}$. Now observe that by Theorem~\ref{modelrep} the $RP$-cofibrant representations are charachterized by the property that for every vertex $j\in \mathcal{Q}_0$ the natural morphism of $R$-modules
\[ \bigoplus\limits_{\stackrel{\alpha\in \mathcal{Q}_1}{\tau(\alpha) = j}}e_{\sigma(\alpha)}(M) \to e_j(M) \]
is injective. To conclude it is sufficient to notice that a morphism $M\to N$ of $\Lambda^{op}$-modules is surjective if and only if for every vertex $j\in \mathcal{Q}_0$ the induced morphism of $R$-modules $e_j(M)\to e_j(N)$ is surjective. The second part of the statement is dual.
\end{proof}
\end{theorem}

\begin{lemma}\label{1gor}
Let $\mathcal{Q}$ be a finite acyclic quiver and let $R$ be a quasi-Frobenius ring. Then the (opposite) path algebra $\Lambda^{op} = R\mathcal{Q}^{op}$ is $1$-Gorenstein.
\begin{proof}
First, we show that $\Lambda^{op}$ has injective-dimension less than $2$ as a left module over itself, i.e. the inequality $\injdim_{\Mod(\Lambda^{op})}\Lambda^{op}\leq 1$ holds.
To begin with we notice that the finite set $\left\{e_j\right\}_{j\in\mathcal{Q}_0}$ of lazy paths in $\Lambda^{op}$ is a collection of primitive orthogonal idempotents such that
\[ 1=\sum\limits_{j\in\mathcal{Q}_0} e_j \]
in $\Lambda^{op}$. It follows that the path algebra $\Lambda^{op}$ itself is isomorphic to the direct sum of a finite number of indecomposable projective $\Lambda^{op}$-modules
\[ \bigoplus_{j\in\mathcal{Q}_0}e_j \Lambda^{op} \]
as a module over $\Lambda^{op}$. Therefore, through the equivalence $\Mod(\Lambda^{op})\simeq \Rep(\mathcal{Q}, R)$, $\Lambda^{op}$ corresponds to a representation $A\in\Rep(\mathcal{Q},R)$ given by the direct sum of some indecomposable projective representations. In particular, $A_j$ is a projective $R$-module for every vertex $j\in Q_0$. Now, consider the morphism $0\xrightarrow{\alpha} A$ in $\Rep(\mathcal{Q},R)$. Clearly, $\alpha$ is vertexwise a stable cofibration of $R$-modules, and then it is a $RI$-stable cofibration by Proposition~\ref{Rstablefib}. Hence, by Proposition~\ref{finiteRTF} we have
\[ \injdim_{\Mod(\Lambda^{op})}\Lambda^{op} = \injdim_{\Rep(\mathcal{Q},R)}\left(A\right) = \injdim_{\Rep(\mathcal{Q},R)}\left(\coker{\alpha}\right) \leq 1 \]
as required.
To prove the inequality $\injdim_{\Mod(\Lambda)}\Lambda \leq 1$ it is sufficient to consider the opposite quiver $\mathcal{Q}^{op}$, and to repeat the above discussion.
\end{proof}
\end{lemma}

\begin{theorem}\label{modelcoincide}
Let $\mathcal{Q}$ be a finite acyclic quiver and let $R$ be a quasi-Frobenius ring. Consider the path algebra $\Lambda = R\mathcal{Q}$. Then the Hovey-projective model structure of Theorem~\ref{Hovey} and the Reedy-projective model structure of Theorem~\ref{modelmodule} coincide. Dually, the Hovey-injective model structure of Theorem~\ref{Hovey} and the Reedy-injective model structure of Theorem~\ref{modelmodule} coincide.
\begin{proof}
It suffices to observe that the class of fibrations and trivial fibrations of the two model structures coincide. In both structures fibrations are surjective morphisms. Since by Lemma~\ref{1gor} $\Lambda^{op}$ is $1$-Gorenstein, it follows that every $\Lambda^{op}$-module of finite projective dimension has projective dimension at most $1$ (see e.g. \cite[Proposition 10.1.15]{EJ11}). Hence by Proposition~\ref{finiteRTF} the trivial fibrations of Theorem~\ref{modelmodule} are surjective morphisms with kernel of finite projective dimension. The second part of the statement is dual.
\end{proof}
\end{theorem}

As an immediate consequence of Theorem~\ref{modelcoincide} we obtain a characterization of Gorenstein-projective and Gorenstein-injective $\Lambda^{op}$-modules.

\begin{corollary}\label{gorcof}
Let $\mathcal{Q}$ be a finite acyclic quiver and let $R$ be a quasi-Frobenius ring. Consider the path algebra $\Lambda = R\mathcal{Q}$. A module $M\in\Mod(\Lambda^{op})$ is Gorenstein-projective if and only if the corresponding $R$-representation is $RP$-cofibrant in $\Rep(\mathcal{Q},R)$. Dually, a module $M\in\Mod(\Lambda^{op})$ is Gorenstein-injective if and only if the corresponding $R$-representation is $RI$-fibrant in $\Rep(\mathcal{Q},R)$.
\begin{proof}
It follows by Theorem~\ref{modelrep}, Theorem~\ref{Hovey}, and Theorem~\ref{modelcoincide}.
\end{proof}
\end{corollary}

This result agrees with the one obtained by Luo and Zhang in \cite{LZ13} since (when $R$ is a finite dimensional self-injective algebra) their definition of \emph{monic representation} precisely coincide with our notion of $RP$-cofibrant representation.
Also, one can then characterize $R$-representations corresponding to Gorenstein-projective (Gorenstein-injective) $\Lambda^{op}$-modules in terms of a lifting property in the category $\Rep(\mathcal{Q},R)$.

\begin{corollary}
Let $\mathcal{Q}$ be a finite acyclic quiver, and let $R$ be a quasi-Frobenius ring. Consider the path algebra $\Lambda=R\mathcal{Q}$ of $\mathcal{Q}$ over $R$, and a representation $M\in\Rep(\mathcal{Q},R)$. Then $M$ corresponds to a Gorenstein-projective $\Lambda^{op}$-module if and only if for every $RP$-stable fibration $p\colon X\to Y$ in $\Rep(\mathcal{Q},R)$ every morphism $q\colon M\to Y$ admits a lifting $h\colon M\to Y$, i.e. the following diagram commutes
\[ \xymatrix{  & X\ar@{->>}[d]_{p}^{\sim} \\
M\ar@{->}[r]_{q}\ar@{.>}[ru]^{h} & Y } \]
in the category $\Rep(\mathcal{Q},R)$.
Dually, $M$ corresponds to a Gorenstein-injective $\Lambda^{op}$-module if and only if for every $RI$-stable cofibration $\iota\colon X\to Y$ in $\Rep(\mathcal{Q},R)$ every morphism $q\colon X\to M$ admits a lifting $h\colon Y\to M$, i.e. the following diagram commutes
\[ \xymatrix{  X\ar@{->}[d]_{\iota}^{\sim} \ar@{->}[r]^{q} & M \\
Y\ar@{.>}[ru]_{h} & } \]
in the category $\Rep(\mathcal{Q},R)$.
\begin{proof}
The statement is an immediate consequence of Corollary~\ref{gorcof}.
\end{proof}
\end{corollary}

\begin{example}
Let $\mathcal{Q}$ be the quiver $0\to 1$, and let $R = \mathbb{C}[\epsilon]$ be the algebra of dual numbers over $\mathbb{C}$, i.e. $R = \mathbb{C}[t]/(t^2)$ where $t$ is a central variable. Since $\mathbb{C}[\epsilon]$ is self-injective we have a model structure on the category $\Rep(\mathcal{Q},\mathbb{C}[\epsilon])$ of representations of $\mathcal{Q}$ over the algebra $\mathbb{C}[\epsilon]$. We now turn our attention to the cofibrant objects. Thanks to Theorem~\ref{modelrep}, a $\mathbb{C}[\epsilon]$-representation $M\in\Rep(\mathcal{Q},\mathbb{C}[\epsilon])$ is $RP$-cofibrant if and only if the morphism of $\mathbb{C}[\epsilon]$-modules $m_{01}\colon M_0\to M_1$ is injective. By Corollary~\ref{gorcof} we have
\[ \Rep(\mathcal{Q},\mathbb{C}[\epsilon])^{RPcof} \simeq \GProj(\mathbb{C}\mathcal{Q}^{op}\otimes_{\mathbb{C}}\mathbb{C}[\epsilon]). \]
Dually, again by Theorem~\ref{modelrep}, a $\mathbb{C}[\epsilon]$-representation $M\in\Rep(\mathcal{Q},\mathbb{C}[\epsilon])$ is $RI$-fibrant if and only if the morphism of $\mathbb{C}[\epsilon]$-modules $m_{01}\colon M_0\to M_1$ is surjective. By Corollary~\ref{gorcof} we have
\[ \Rep(\mathcal{Q},\mathbb{C}[\epsilon])^{RIfib} \simeq \GInj(\mathbb{C}\mathcal{Q}^{op}\otimes_{\mathbb{C}}\mathbb{C}[\epsilon]). \]
\end{example}

Thanks to Corollary~\ref{gorcof} one can also explicitly characterize Reedy stable equivalences in $\Rep(\mathcal{Q},R)$. This in fact will be useful in Section~\ref{homotopy}.

\begin{proposition}\label{Reedyeq}
Let $\mathcal{Q}$ be a finite acyclic quiver, and let $R$ be a quasi-Frobenius ring. A morphism $f\colon M\to N$ in $\Rep(\mathcal{Q},R)$ is a Reedy stable equivalence if and only if there exists a commutative diagram
\[ \xymatrix{  G_PM \ar@{->}[r]^{\iota} \ar@{->}[d] & G_PM\oplus T \ar@{->}^{\cong}[r] & G_PN\oplus S \ar@{->}[r]^{\pi} & G_PN \ar@{->}[d] \\
M\ar@{->}^{f}[rrr] & & & N } \]
for some projective representations $S,T\in\Rep(\mathcal{Q},R)$, where $G_PM\to M$ and $G_PN\to N$ are $RP$-cofibrant replacements for $M$ and $N$ respectively, $\iota$ is the natural inclusion and $\pi$ is the natural projection.

Dually, a morphism $f\colon M\to N$ in $\Rep(\mathcal{Q},R)$ is a Reedy stable equivalence if and only if there exists a commutative diagram
\[ \xymatrix{  M\ar@{->}^{f}[rrr] \ar@{->}[d] & & & N \ar@{->}[d] \\
G_IM \ar@{->}[r]^{\iota} & G_IM\oplus J \ar@{->}^{\cong}[r] & G_IN\oplus K \ar@{->}[r]^{\pi} & G_IN  } \]
for some injective representations $J,K\in\Rep(\mathcal{Q},R)$, where $M\to G_IM$ and $N \to G_IN$ are $RI$-fibrant replacements for $M$ and $N$ respectively, $\iota$ is the natural inclusion and $\pi$ is the natural projection.
\begin{proof}
Consider a Reedy stable equivalence $f\colon M\to N$ in $\Rep(\mathcal{Q},R)$. Applying the $RP$-cofibrant replacement functor we obtain a commutative square
\[ \xymatrix{  G_PM \ar@{->}[r]^{G_Pf} \ar@{->}[d] & G_PN \ar@{->}[d] \\
M\ar@{->}^{f}[r] & N. } \]
Now take a factorization of $G_Pf$ as a $RP$-stable cofibration followed by a $RP$-fibration:
\[ G_Pf\colon G_PM \xrightarrow{\alpha} D\xrightarrow{\beta} G_PN. \]
By the two-out-of-three axiom $\beta$ is a $RP$-stable fibration. Now, $\alpha$ fits into an exact sequence $0\to G_PM\xrightarrow{\alpha} D \to \coker(\alpha)\to 0$, and this sequence splits since $G_PM$ is $RP$-fibrant. Then, it is easy to see that the left lifting property of $\alpha$ with respect to $RP$-fibrations is equivalent to the projectivity of $\coker(\alpha)$ in the category $\Rep(\mathcal{Q},R)$. Also, since $G_PN$ is $RP$-cofibrant $\beta$ splits too. This gives a morphism $D\xrightarrow{r} \ker(\beta)$ such that the composition $\ker(\beta)\hookrightarrow D\xrightarrow{r} \ker(\beta)$ is the identity morphism on $\ker(\beta)$, so that $\ker(\beta)$ is a retract of $D$ and hence $RP$-cofibrant. Being also a trivial object in $\Rep(\mathcal{Q},R)$, $\ker(\beta)$ is projective by Theorem~\ref{Hovey}. The statement follows by taking $T=\coker(\alpha)$ and $S=\ker(\beta)$. The second part of the statement is dual.
\end{proof}
\end{proposition}

\section{The stable category of quiver representations}\label{homotopy}

The aim of this section is to investigate the main properties of the homotopy category of $\Rep(\mathcal{Q},R)$. As we will see, it satisfies three different universal properties. In particular it is an \emph{algebraic} category, meaning that it is triangle equivalent to the stable category of a Frobenius category, namely to $\underline{\GProj}(\Lambda^{op})$ and $\underline{\GInj}(\Lambda^{op})$. We begin by recalling definitions and notations.

A \emph{Frobenius} category is a Quillen exact category which has enough injectives and enough projectives, and where the class of projectives coincides with the class of injectives.
For instance, given a quasi-Frobenius ring, the category $\Mod(R)$ is a Frobenius category where the Quillen exact structure is given by the short exact sequences. Another interesting example is given by the following result, which actually holds in a more general setting; for details and relations with Geometry we refer to~\cite{IKWY}.

\begin{proposition}
Let $\mathcal{Q}$ be a finite acyclic quiver and let $R$ be a quasi-Frobenius ring. Consider the path algebra $\Lambda = R\mathcal{Q}$. Then the full subcategory $\GProj(\Lambda^{op})\subseteq\Mod(\Lambda^{op})$ of Gorenstein-projective $\Lambda^{op}$-modules is a Frobenius category. Dually, the full subcategory $\GInj(\Lambda^{op})\subseteq\Mod(\Lambda^{op})$ of Gorenstein-injective $\Lambda^{op}$-modules is a Frobenius category.
\end{proposition}

It is well-known that the stable category of a Frobenius category is triangulated, see \cite{HZ12}. Given a Frobenius category $\mathbf{F}$, its stable category $\underline{\mathbf{F}}$ is defined as follows. The objects are the same as $\mathbf{F}$, while given $A, B\in \underline{\mathbf{F}}$ one defines $\Hom_{\underline{\mathbf{F}}}(A,B)=\frac{\Hom_{\mathbf{F}}(A,B)}{J_{A,B}}$, where $J_{A,B}$ is the ideal generated by all those morphisms $A\to B$ which factor through an injective/projective object in $\mathbf{F}$. There exists an obvious functor $\gamma\colon \mathbf{F}\to\underline{\mathbf{F}}$ that is the identity on objects. In fact, the pair $(\underline{\mathbf{F}},\gamma)$ is initial in the category of categories under $\mathbf{F}$ annihilating injective/projective objects. In other terms, $\underline{\mathbf{F}}$ satisfies the following universal property: for every functor $\alpha\colon\mathbf{F}\to \mathbf{C}$ such that $\alpha(P)\cong\alpha(0_{\mathbf{F}})$ for every projective $P\in\mathbf{F}$, there exists a functor $\widetilde{\alpha}\colon\underline{\mathbf{F}}\to\mathbf{C}$ such that $\alpha = \widetilde{\alpha}\circ\gamma$.

The shift functor $T\colon \underline{\mathbf{F}}\to \underline{\mathbf{F}}$ is defined as follows. Take an object $A\in\underline{\mathbf{F}}$ and consider a conflation
$A\to I_A\to T_A$
where $I_A$ is an injective object in $\mathbf{F}$.
Our interest will be only in categories where conflations are precisely short exact sequences. Then, $T(A)=T_A$ is well-defined in $\underline{\mathbf{F}}$ by the Schanuel's Lemma (for details we refer again to \cite{HZ12}).

Now recall that in the homotopy category of any (pointed) model category $\mathbf{M}$ the suspension functor $\Sigma$ and the loop functor $\Omega$ are defined. More precisely, given an object $A\in\mathbf{M}$, we first take the cofibrant replacement $C_A\to A$, then we factor the morphism $C_A\to 0$ as a cofibration followed by a trivial fibration $C_A\xrightarrow{\iota} D_A\to 0$, so that we can define $\Sigma A=\coker(\iota)$. Dually, given an object $A\in\mathbf{M}$, we first take the fibrant replacement $A\to F_A$, then factor the morphism $0\to F_A$ as a trivial cofibration followed by a fibration $0\to E_A\xrightarrow{\pi} F_A$, and define $\Omega A = \ker(\pi)$. Notice that up to weak equivalences all the choices we made do not really matter.

Now observe that the identity functor $\Rep(\mathcal{Q},R)\to\Rep(\mathcal{Q},R)$ is a Quillen equivalence from the $RP$-model structure to the $RI$-model structure of Theorem~\ref{modelrep}, so that the homotopy categories of these two model structures are equivalent. In fact it is immediate to check that they are isomorphic. From now on we denote by $\Ho(\Rep(\mathcal{Q},R))$ the homotopy category with respect to both Reedy model structures of Theorem~\ref{modelrep}. Similarly, in the following we denote by $\Ho(\Mod(G))$ the homotopy category of modules over a Gorenstein ring with respect to both Hovey model structures of Theorem~\ref{Hovey}.
Since the suspension functor and the loop functor are preserved under Quillen equivalences, we can construct $\Sigma$ and $\Omega$ in either the $RP$-model structure or $RI$-model structure. Hence to construct the suspension of a given $R$-representation $M\in\Rep(\mathcal{Q},R)$ we begin by observing that it is $RI$-cofibrant (since everything is so) and then we consider an exact sequence
\[ 0\to M\xrightarrow{\iota} I_M\to \coker(\iota)\to 0 \]
where $I_M$ is an injective $R$-representation, and $\iota$ is a $RI$-cofibration (i.e. vertexwise injective). We define $\Sigma M =\coker(\iota)$. Dually, using the $RP$-model structure on $\Rep(\mathcal{Q},R)$ we can construct the loop $\Omega M$ as the $R$-representation fitting in a short exact sequence
\[ 0\to \Omega M\to P_M\to M \to 0 \]
with $P_M$ projective in $\Rep(\mathcal{Q},R)$.

Our aim is now to show that $\Ho(\Rep(\mathcal{Q},R))$ is triangle equivalent to the stable categories $\underline{\GProj}(\Lambda^{op})$ and $\underline{\GInj}(\Lambda^{op})$; as a consequence it is an algebraic category which satisfies three different universal properties. As we will see, the same holds for $\underline{\GProj}(\Lambda^{op})$ and for $\underline{\GInj}(\Lambda^{op})$. We first need to recall the following result.

\begin{theorem}[Hovey, \cite{Hov02}]
Let $G$ be a $1$-Gorenstein ring. Then $\Sigma$ and $\Omega$ are inverse equivalences of $\Ho(\Mod(G))$, therefore $\Ho(\Mod(G))$ is triangulated and $\Sigma$ is the shift functor. Moreover:
\[ \Hom_{\Ho(\Mod(G))}(\Omega M,N) \cong \Hom_{\Ho(\Mod(G))}(M,\Sigma N) \cong \Ext^1(G_PM, N)\cong\Ext^1(M,G_IN) \]
where $G_PM$ is a $RP$-cofibrant replacement for $M$ while $G_IN$ is a $RI$-fibrant replacement for N.
\end{theorem}

\begin{theorem}\label{equivalences}
Let $\mathcal{Q}$ be a finite acyclic quiver and let $R$ be a quasi-Frobenius ring. Consider the path algebra $\Lambda = R\mathcal{Q}$. Then there exist equivalences of triangulated categories
\[ \underline{\GProj}(\Lambda^{op})\simeq \Ho(\Rep(\mathcal{Q},R)) \simeq \underline{\GInj}(\Lambda^{op}). \]
In particular, $\Ho(\Rep(\mathcal{Q},R))$ is an algebraic category.
\begin{proof}
Consider the composite functor
\[ \delta\colon \GProj(\Lambda^{op})\hookrightarrow \Mod(\Lambda^{op}) \xrightarrow{\simeq} \Rep(\mathcal{Q},R) \xrightarrow{\gamma_R} \Ho(\Rep(\mathcal{Q},R)) \]
where $\Rep(\mathcal{Q},R) \xrightarrow{\gamma_R} \Ho(\Rep(\mathcal{Q},R))$ is the projection on the homotopy category with respect to the $RP$-model structure.
Clearly $\delta$ sends all the projective $\Lambda^{op}$-modules to zero in $\Ho(\Rep(\mathcal{Q},R))$ and then there exists a lifting functor $\widetilde{\delta}\colon \underline{\GProj}(\Lambda^{op})\to \Ho(\Rep(\mathcal{Q},R))$ such that $\delta = \widetilde{\delta}\circ\gamma_{\Lambda}$ where $\gamma_{\Lambda}\colon\GProj(\Lambda^{op})\to \underline{\GProj}(\Lambda^{op})$ is the projection on the stable category. One can easily see that $\widetilde{\delta}$ is the identity on objects and that it is a dense functor since every $R$-representation is naturally isomorphic to its $RP$-cofibrant replacement (which is Gorenstein-projective). Anyway it is possible to explicitly exhibit the quasi-inverse functor of $\widetilde{\delta}$. Indeed, consider the composite functor
\[ \omega\colon\Rep(\mathcal{Q},R) \xrightarrow{G_P} \GProj(\Lambda^{op}) \xrightarrow{\gamma_{\Lambda}} \underline{\GProj}(\Lambda^{op}) \]
where $\Rep(\mathcal{Q},R) \xrightarrow{G_P} \GProj(\Lambda^{op})$ is the $RP$-cofibrant replacement functor, while the functor $\GProj(\Lambda^{op}) \xrightarrow{\gamma_{\Lambda}} \underline{\GProj}(\Lambda^{op})$ is simply the projection on the stable category.
By Proposition~\ref{Reedyeq} it is immediate to see that every $RP$-stable equivalence in $\Rep(\mathcal{Q},R)$ is sent to an isomorphism in $\underline{\GProj}(\Lambda^{op})$. Hence, by the universal property of $\Ho(\Rep(\mathcal{Q},R))$ there exists a lifting functor $\widetilde{\omega}\colon \Ho(\Rep(\mathcal{Q},R)) \to \underline{\GProj}(\Lambda^{op})$ such that $\omega = \widetilde{\omega} \circ \gamma_R$. It is now straightforward to check that $\widetilde{\omega}$ is the required quasi-inverse for $\widetilde{\delta}$.

The second equivalence of categories can be proved dually. To conclude, notice that the suspension functor $\Sigma$ in $\Ho(\Rep(\mathcal{Q},R))$ has precisely the same description as the shift functor in the stable category of the Frobenius categories $\GProj(\Lambda^{op})$ and $\GInj(\Lambda^{op})$, so that it is preserved by $\widetilde{\delta}$.
\end{proof}
\end{theorem}

Theorem~\ref{equivalences} ``allows" the homotopy category $\Ho(\Rep(\mathcal{Q},R))$ to be called the \emph{stable category of quiver representations over} $R$.

\begin{remark}
One of the possible ``models" (up to equivalences of categories) for the homotopy category $\Ho(\mathbf{M})$ of a model category $\mathbf{M}$ can be constructed as follows. Consider the full subcategory $\mathbf{M}^{cf}\subseteq\mathbf{M}$ whose objects are those of $\mathbf{M}$ that are both fibrant and cofibrant. Given two objects $A,B\in\mathbf{M}^{cf}$ we consider the homotopy relation $\sim_h$ on $\Hom_{\mathbf{M}^{cf}}(A,B)$, see \cite[Definition 1.2.4]{Hov99}.
Then one defines the category $\underline{\mathbf{M}}^{cf}$ as follows:
\begin{description}
\item[1] $\Ob(\underline{\mathbf{M}}^{cf}) = \Ob(\mathbf{M}^{cf})$,
\item[2] $\Hom_{\underline{\mathbf{M}}^{cf}}(A,B) = \frac{\Hom_{\mathbf{M}^{cf}}(A,B)}{\sim_h}$ for any $A,B\in\underline{\mathbf{M}}^{cf}$.
\end{description}
One of the main results in Model Category Theory, the so-called \emph{fundamental theorem} of model categories, states that the inclusion $\mathbf{M}^{cf}\hookrightarrow \mathbf{M}$ induces an equivalence of categories $\underline{\mathbf{M}}^{cf}\simeq\Ho(\mathbf{M})$, see \cite[Theorem 1.2.10]{Hov99}.

Notice that in the $RP$-model structure of Theorem~\ref{modelrep} every $R$-representation is $RP$-fibrant and then the subcategory of $RP$-fibrant-cofibrant objects is $\Rep(\mathcal{Q},R)^{RPcof}$ which corresponds to Gorenstein-projective modules. Looking at the proof of Theorem~\ref{equivalences} one sees that in fact the equivalence of categories $\Rep(\mathcal{Q},R)\xrightarrow{\simeq}\Mod(\Lambda^{op})$ is, in particular, a Quillen equivalence and that the homotopy relation between Gorenstein-projective modules is precisely the relation induced by declaring two morphisms $f,g\in\Hom_{\GProj(\Lambda^{op})}(A,B)$ stably equivalent if $(f-g)$ factors through a projective module in $\GProj(\Lambda^{op})$.
\end{remark}

\begin{corollary}
Let $\mathcal{Q}$ be a finite acyclic quiver and let $R$ be a quasi-Frobenius ring. Consider the path algebra $\Lambda = R\mathcal{Q}$. Then the homotopy category $\Ho(\Rep(\mathcal{Q},R))$ satisfies the following universal properties.
\begin{description}
\item[1] Consider the composite functor
\[ \gamma_P\colon\GProj(\Lambda^{op})\hookrightarrow\Mod(\Lambda^{op}) \xrightarrow{\simeq} \Rep(\mathcal{Q},R)\xrightarrow{\gamma_R}\Ho(\Rep(\mathcal{Q},R)) \]
where $\Rep(\mathcal{Q},R) \xrightarrow{\gamma_R} \Ho(\Rep(\mathcal{Q},R))$ is the projection on the homotopy category with respect to the $RP$-model structure.
Given a functor $\alpha\colon\GProj(\Lambda^{op})\to \mathbf{C}$ such that for every projective module $Q\in\GProj(\Lambda^{op})$ there is an isomorphism $\alpha(Q)\cong\alpha(0)$ in $\mathbf{C}$, there exists a lifting $\widetilde{\alpha}\colon\Ho(\Rep(\mathcal{Q},R))\to\mathbf{C}$ such that the diagram
\[ \xymatrix{ \GProj(\Lambda^{op})\ar@{->}[r]^{\alpha} \ar@{->}[d]_{\gamma_P} & \mathbf{C} \\
\Ho(\Rep(\mathcal{Q},R)) \ar@{.>}_{\widetilde{\alpha}}[ur] & } \]
commutes.
\item[2] Consider the composite functor
\[ \gamma_I\colon\GInj(\Lambda^{op})\hookrightarrow\Mod(\Lambda^{op}) \xrightarrow{\simeq} \Rep(\mathcal{Q},R)\xrightarrow{\gamma_R'}\Ho(\Rep(\mathcal{Q},R)) \]
where $\Rep(\mathcal{Q},R) \xrightarrow{\gamma_R'} \Ho(\Rep(\mathcal{Q},R))$ is the projection on the homotopy category with respect to the $RI$-model structure.
Given a functor $\alpha\colon\GInj(\Lambda^{op})\to \mathbf{C}$ such that for every injective module $J\in\GInj(\Lambda^{op})$ there is an isomorphism $\alpha(J)\cong\alpha(0)$ in $\mathbf{C}$, there exists a lifting $\widetilde{\alpha}\colon\Ho(\Rep(\mathcal{Q},R))\to\mathbf{C}$ such that the diagram
\[ \xymatrix{ \GInj(\Lambda^{op})\ar@{->}[r]^{\alpha} \ar@{->}[d]_{\gamma_I} & \mathbf{C} \\
\Ho(\Rep(\mathcal{Q},R)) \ar@{.>}_{\widetilde{\alpha}}[ur] & } \]
commutes.
\end{description}
\begin{proof}
It is an immediate consequence of Theorem~\ref{equivalences}.
\end{proof}
\end{corollary}

Notice that, thanks to Theorem~\ref{equivalences}, $\underline{\GProj}(\Lambda^{op})$ inherits the universal properties of $\underline{\GInj}(\Lambda^{op})$ and $\Ho(\Rep(\mathcal{Q},R))$.

\pagestyle{plain}

\end{document}